\documentclass[11pt,a4paper]{amsart}

\usepackage{latexsym}
\usepackage{amsmath}
\usepackage{amsthm}
\usepackage{amssymb}
\usepackage{vmargin}
\usepackage{amscd}
\usepackage{stmaryrd}
\usepackage{euscript}
\usepackage{mathrsfs}
\usepackage{amscd}
\usepackage[all]{xy}
\usepackage{xr}

\DeclareMathAlphabet{\mathpzc}{OT1}{pzc}{m}{it}

\newtheorem{theorem}{Theorem}[section]
\newtheorem{proposition}[theorem]{Proposition}
\newtheorem{corollary}[theorem]{Corollary}

\newtheorem{lemma}[theorem]{Lemma}
\newtheorem*{theorem*}{Theorem}
\newtheorem*{proposition*}{Proposition}
\newtheorem*{corollary*}{Corollary}
\newtheorem*{lemma*}{Lemma}
\newtheorem*{conjecture*}{Conjecture}

\theoremstyle{definition}
\newtheorem{definition}[theorem]{Definition}

\newtheorem*{definition*}{Definition}

\theoremstyle{remark}
\newtheorem{example}[theorem]{Example}

\newtheorem{remark}[theorem]{Remark}
\newtheorem{remarks}[theorem]{Remarks}

\newtheorem*{example*}{Example}
\newtheorem*{examples*}{Examples}
\newtheorem*{remark*}{Remark}
\newtheorem*{remarks*}{Remarks}
\newtheorem*{exercise*}{Exercise}

\newcommand\da{\!\downarrow\!}

\newcommand\la{\leftarrow}

\newcommand\id{\mathrm{id}}

\newcommand\ten{\otimes}

\newcommand\CC{\mathrm{C}}
\newcommand\DD{\mathrm{D}}

\renewcommand\H{\mathrm{H}}

\newcommand\N{\mathbb{N}}
\newcommand\Z{\mathbb{Z}}
\newcommand\Q{\mathbb{Q}}

\newcommand\bG{\mathbb{G}}

\newcommand\bL{\mathbb{L}}

\newcommand\bS{\mathbb{S}}

\newcommand\C{\mathcal{C}}

\newcommand\cA{\mathcal{A}}
\newcommand\cB{\mathcal{B}}

\newcommand\cD{\mathcal{D}}
\newcommand\cE{\mathcal{E}}
\newcommand\cF{\mathcal{F}}
\newcommand\cG{\mathcal{G}}

\newcommand\cL{\mathcal{L}}
\newcommand\cM{\mathcal{M}}
\newcommand\cN{\mathcal{N}}

\newcommand\cT{\mathcal{T}}

\newcommand\cW{\mathcal{W}}

\newcommand\D{\mathcal{D}}

\newcommand\F{\mathscr{F}}

\renewcommand\O{\mathscr{O}}

\newcommand\sC{\mathscr{C}}

\newcommand\sE{\mathscr{E}}
\newcommand\sF{\mathscr{F}}

\newcommand\sH{\mathscr{H}}

\newcommand\sU{\mathscr{U}}

\newcommand\fM{\mathfrak{M}}

\newcommand\fX{\mathfrak{X}}
\newcommand\fY{\mathfrak{Y}}
\newcommand\fZ{\mathfrak{Z}}

\renewcommand\L{\Lambda}

\newcommand\m{\mathfrak{m}}

\newcommand\cHom{\mathcal{H}\!\mathit{om}}

\newcommand\sEnd{\mathscr{E}\!\mathit{nd}}

\newcommand\Ho{\mathrm{Ho}}
\newcommand\Ring{\mathrm{Ring}}
\newcommand\Alg{\mathrm{Alg}}
\newcommand\ALG{\mathrm{ALG}}
\newcommand\FALG{\mathrm{FALG}}
\newcommand\CART{\mathrm{CART}}
\newcommand\HYP{\mathrm{HYP}}
\newcommand\hyp{\mathrm{hyp}}
\newcommand\Mod{\mathrm{Mod}}

\newcommand\Hom{\mathrm{Hom}}

\newcommand\HHom{\underline{\mathrm{Hom}}}
\newcommand\EEnd{\underline{\mathrm{End}}}

\newcommand\Ext{\mathrm{Ext}}
\newcommand\EExt{\mathbb{E}\mathrm{xt}}
\newcommand\End{\mathrm{End}}

\newcommand\Aut{\mathrm{Aut}}

\newcommand\cone{\mathrm{cone}}

\newcommand\coker{\mathrm{coker\,}}

\newcommand\im{\mathrm{Im\,}}

\newcommand\Ob{\mathrm{Ob}\,}

\newcommand\Ab{\mathrm{Ab}}

\newcommand\Gp{\mathrm{Gp}}

\newcommand\Stack{\mathrm{Stack}}

\newcommand\Spec{\mathrm{Spec}\,}
\newcommand\Dec{\mathrm{Dec}\,}
\newcommand\DEC{\mathrm{DEC}\,}

\newcommand\Set{\mathrm{Set}}

\newcommand\Cat{\mathrm{Cat}}

\newcommand\Aff{\mathrm{Aff}}

\newcommand\Sp{\mathrm{Sp}}
\newcommand\Sch{\mathrm{Sch}}

\newcommand\ad{\mathrm{ad}}

\newcommand\Lim{\varprojlim}
\newcommand\LLim{\varinjlim}

\newcommand\ho{\mathrm{ho}\!}
\newcommand\into{\hookrightarrow}
\newcommand\onto{\twoheadrightarrow}
\newcommand\abuts{\implies}
\newcommand\xra{\xrightarrow}
\newcommand\xla{\xleftarrow}

\newcommand\bt{\bullet}
\newcommand\by{\times}

\newcommand\Symm{\mathrm{Symm}}
\newcommand\SL{\mathrm{SL}}
\newcommand\GL{\mathrm{GL}}

\newcommand\cart{\mathrm{cart}}

\newcommand\diag{\mathrm{diag}\,}

\newcommand\ind{\mathrm{ind}}

\newcommand\pd{\partial}

\newcommand\Tor{\mathrm{Tor}}

\newcommand\st{\rm{st}}

\newcommand\gp{\mathrm{Gp}}
\newcommand\gpd{\mathrm{Gpd}}

\newcommand\op{\mathrm{opp}}

\newcommand\co{\colon\thinspace}

\newcommand\oR{\mathbf{R}}

\newcommand\oL{\mathbf{L}}

\newcommand\uleft\underleftarrow
\newcommand\uline\underline
\newcommand\uright\underrightarrow

\newcommand\HOM{\mathrm{HOM}}

\begin{document}

\begin{abstract}
We describe derived moduli functors for a range of problems involving schemes and quasi-coherent sheaves, and give cohomological conditions for them to be representable by derived geometric $n$-stacks. Examples of problems represented by derived geometric $1$-stacks are derived moduli of polarised projective varieties, derived moduli of vector bundles, and derived moduli of 
abelian varieties.
\end{abstract}

\title{Derived moduli of schemes and sheaves}

\author{J.P.Pridham}
\thanks{The author was supported during this research by  the Engineering and Physical Sciences Research Council [grant number  EP/F043570/1].}
\maketitle

\section*{Introduction}

The driving philosophy behind derived algebraic geometry is that for every moduli problem, there should be an associated derived moduli problem. Thanks to  the framework of \cite{hag2} and \cite{lurie}, there is often an obvious candidate for such a functor for derived moduli of schemes:

 Given a simplicial algebra $A$, we first form the $\infty$-category $\C(A)$ of derived geometric stacks $\fX$   homotopy-flat over $\Spec A$, with the additional property that
\[
\fX\by^h_{\Spec A} \Spec \pi_0A
\] 
is a scheme. We can impose additional  restrictions (smoothness, properness, dimension) on $\fX$ by requiring that they apply to $\fX\by^h_{\Spec A} \Spec \pi_0A$. We then consider  the $\infty$-category $\cG(A)\subset \C(A)$ in which the only morphisms are weak equivalences. 
Finally, we define the derived moduli functor $\fM(A)$ to be the nerve of the $\infty$-category $\cG(A)$. There are variants of this construction for subschemes and so on.

Similarly, for derived moduli of sheaves or complexes on a scheme $X$, we can first set $\C(A)$ to be the $\infty$-category of quasi-coherent complexes $\sF$ (in the sense of \cite{hag2}) on $X\by \Spec A$,  such that the complex $\sF\ten^{\oL}_{A}\pi_0A$ on $X \by \Spec (\pi_0A)$ satisfies any additional constraints (bounded, fixed Euler characteristic, concentrated in degree $0$, free of rank $n$). We then form $\cG(A) \subset  \C(A)$ as before, and take the nerve.

We can even combine these approaches to consider moduli of polarised projective schemes. $\C(A)$ then consists of  pairs $(\fX,\sF)$, for $\fX$ a derived geometric stack    homotopy-flat over $\Spec A$ and $\sF$ a quasi-coherent complex on $\fX$, with 
\[
(\fX\by^h_{\Spec A} \Spec \pi_0A, \sF\ten^{\oL}_{A}\pi_0A)
\]
a projective scheme and an ample line bundle.

Lurie's Representability Theorem provides explicit criteria to check whether such a functor is representable by a geometric derived stack. However, the framework of \cite{hag2} and \cite{lurie} does not lend itself easily to calculations, making it difficult to verify that these functors satisfy the criteria of the representability theorem. Even to show that the cotangent complex governs first-order deformations of $n$-geometric stacks, the more concrete approach of Artin hypergroupoids is needed (as in \cite{stacks2} \S \ref{stacks-cotsn}). 

Several examples of representable derived moduli functors have already been established, however. Both \cite{hag2} and \cite{lurie} construct  $\Hom$-stacks, while \cite{lurie} 
addresses derived moduli of stable curves, as well as derived Picard stacks and derived Hilbert schemes.  \cite{hag2} deals with derived moduli of local systems and of finite algebras for any operad, and \cite{toenseattle} describes derived moduli of stable maps.

The purpose of this paper is to establish representability more generally for problems involving schemes and quasi-coherent sheaves,  building on the explicit characterisations of \cite{stacks2}. In simplistic terms, we show that if the associated underived moduli problem is representable, then the derived moduli problem will be representable provided only that certain  cohomology groups satisfy mild finiteness conditions. The main results are Theorem \ref{representdaffine} for  moduli of derived geometric stacks over a fixed base, Theorem \ref{augrepresentdaffine} which has the additional datum of an augmentation, and Theorem \ref{representdmod} for moduli of quasi-coherent complexes.

In contrast to the indirect approach of satisfying a representability theorem, \cite{Hilb} and \cite{Quot} construct explicit derived Hilbert and Quot schemes as dg-schemes with the necessary  properties, but there  no universal family is given, so the derived moduli spaces lack functorial interpretations. 
In \cite{dmc}, we  use the results of this paper to compare these approaches, thereby giving explicit presentations for the derived moduli spaces constructed here. 

The structure of the paper is as follows.

In Section \ref{simpcat} we review the framework of simplicial categories, giving several properties of the nerve construction from simplicial categories to simplicial sets. Section \ref{modsimpcat} then recalls a variant (Theorem \ref{lurierep3}) of Lurie's Representability Theorem, and shows how the nerve construction transfers good properties from simplicial category-valued functors to simplicial set-valued functors. These results are then applied in the last two sections to a wide range of moduli problems.  

The main results of Section \ref{modschsn} are Corollaries \ref{representaffine} and \ref{representrelaffine}, which deal with derived moduli of affine schemes, and   Theorems \ref{representdaffine} and \ref{augrepresentdaffine}, addressing derived moduli of (augmented) geometric $n$-stacks. These are applied to 
moduli of finite schemes (Example \ref{fineg}),  of schemes (Example \ref{modscheg}), of torsors (Example \ref{Gtorsors}), of polarised schemes (Example \ref{modpolsch}), of 
abelian schemes (Example \ref{abvar}) 
and derived $\Hom$-stacks (Example \ref{Mor}).

Section \ref{modshfsn} deals, in   Theorem \ref{representdmod} and Example \ref{modqcoheg}, with moduli of quasi-coherent complexes.   
Lemma \ref{cfbundles} then checks that the paper's two possible approaches to derived moduli of vector bundles (via $\GL_n$-torsors and via quasi-coherent complexes) give equivalent results, which was already known in the framework of \cite{hag2}.

\tableofcontents


\section{Simplicial categories}\label{simpcat}

Let $\bS$ denote the category of simplicial sets, and $s\bS$ the category of bisimplicial sets.

\begin{definition}
 Let $s\Cat$ be the category of small simplicially enriched  categories, which we will refer to as simplicial categories. Explicitly, an object $\C \in s\Cat$ consists of a set $\Ob \C$ of objects, together with $\HHom_{\C}(x,y) \in \bS$ for all $x,y \in \Ob \bS$, equipped with an associative composition law and identities. 
\end{definition}

\begin{definition}\label{pi0C}
Given a simplicial category $\C$,  the category $\pi_0\C$ is defined to have the same objects as $\C$, with morphisms 
$$
\Hom_{\pi_0\C}(x,y)=\pi_0\HHom_{\C}(x,y). 
$$
\end{definition}

\begin{definition}
Recall from \cite{bergner} Theorem 1.1 that  a morphism  $f:\C \to \D$ in $s\Cat$ is said to be  a  weak equivalence (a.k.a. an $\infty$-equivalence) whenever
\begin{enumerate} 
\item[(W1)] for any objects $a_1$ and $a_2$ in $\C$, the map
$\HHom_{\C}(a_1, a_2)\to \HHom_{\cD}(fa_1, fa_2)$ 
is a weak equivalence of simplicial sets;
\item[(W2)] the induced functor $\pi_0f : \pi_0\C \to \pi_0\cD$ is an equivalence of
categories.
\end{enumerate}
\end{definition}

\begin{definition}
Given a simplicial category $\C$, a morphism in $\C_0$ is said to be a homotopy equivalence if it becomes an isomorphism in $\pi_0\C$. 
\end{definition}

\begin{definition}\label{scatfibdef}
Recall from \cite{bergner} Theorem 1.1 that  a morphism $f:\C \to \cD$ in $s\Cat$ is said to be a fibration whenever
\begin{enumerate} 
\item[(F1)] for any objects $a_1$ and $a_2$ in $\C$, the map
$\HHom_{\C}(a_1, a_2)\to \HHom_{\cD}(fa_1, fa_2)$ 
is a fibration of simplicial sets;
\item[(F2)] for any objects $a_1 \in\C$, $b \in \cD$, and homotopy equivalence $e :
fa_1 \to b$ in $\cD$, there is an object $a_2 \in \C$ and a homotopy equivalence
$d : a_1 \to a_2$ in $\C$ such that $fd = e$.
\end{enumerate}
\end{definition}

\begin{definition}
Given $\C \in s\Cat$, define 
$\bar{W}\C:= \bar{W}B\C$, where $\bar{W}:s\bS \to \bS$ is the right adjoint to  Illusie's total $\Dec$ functor given by $\DEC(X)_{mn}= X_{m+n+1}$, and $B:\Cat \to \bS$ is the nerve. Explicitly
\[
(\bar{W}\C)_n = \coprod_{\uline{x}\in (\Ob\C)^{n+1}}\C_{n-1}(x_{n-1}, x_{n}) \by\C_{n-2}(x_{n-2}, x_{n-1})\by \ldots \by \C_0(x_0, x_1),
\]
with operations
\begin{eqnarray*}
\pd_i( g_{n-1},\ldots ,g_0)&=& 
\left\{ \begin{matrix} 
(g_{n-2},\ldots, g_0) & i=0,\\
(\pd_{i-1}g_{n-1},\ldots, \pd_1g_{n-i+1}, (\pd_0g_{n-i})g_{n-i-1}, g_{n-i-2}, \ldots, g_0) & 0<i<n,\\ 
(\pd_{n-1}g_{n-1}, \ldots, \pd_1g_1) & i=n, 
\end{matrix} \right.\\
\sigma_i( g_{n-1}, \ldots , g_0)&=& ( \sigma_{i-1}g_{n-1},\ldots, \sigma_0g_{n-i}, \id_{x_{n-i}} ,  g_{n-i-1}, \ldots, g_0),
\end{eqnarray*}
\end{definition}

\begin{remark}
Another functor from bisimplicial to simplicial sets is the diagonal functor $\diag(Y)_n= Y_{nn}$.
In \cite{CRdiag},  it is established that the canonical natural transformation
\[
\diag Y \to \bar{W}Y
\]
is a weak equivalence for all $Y$. When $Y=B\C$, this transformation is given by
\[
(y_n, h_{n-1},y_{n-1},\ldots, h_0 ,y_0)\mapsto (y_n, \pd_0h_{n-1}, \pd_0y_{n-1}, \ldots, (\pd_0)^nh_0, (\pd_0)^ny_0).
\]
This observation means that throughout this paper, $\bar{W}$ can be replaced with $\diag$, which is simpler to define. Our reasons for preferring $\bar{W}$ are that it produces much smaller objects, 
as exploited in \cite{pathgpd}. For instance, if $Y$ is a bisimplicial abelian group, then the Dold--Kan normalisation  $N\bar{W}Y$ (see Definition \ref{normdef}) is just the total complex of the Dold--Kan binormalisation $NY$.
\end{remark}

The following results seem to be folklore, but lacking a suitable reference, we prove them afresh.

\begin{proposition}\label{fibscat}
If $f: \C\to \cD$ is a fibration of simplicial categories, with $\pi_0\C$ and $\pi_0\cD$  being groupoids, and $\cB \to \cD$ is an arbitrary morphism in $s\Cat$, then the homotopy fibre product
\[
(\bar{W} \C)\by^h_{(\bar{W}\cD)}(\bar{W}\cB) 
\]
is weakly equivalent to $\bar{W}(\C\by_{\cD}\cB)$.
\end{proposition}
\begin{proof}
Construct a category $\pi_n\C$ with objects $\Ob \C$ and morphisms $\coprod_{f \in \C_0(x,y)}\pi_n(\C(x,y), f)$ for $x, y \in \Ob \C$. If we form the bisimplicial set $B\C$, then $\pi_nB\C= B\pi_n\C$ and we wish to verify the $\pi_*$-Kan condition (\cite{sht} \S IV.4) --- that
\[
B\pi_n\C \to B\C_0
\] 
is a fibration. 

If a functor $F:\cA \to \cA'$ of categories is an isomorphism on objects, then $B\cA \to B\cA'$ is a fibration if and only if 
\begin{enumerate}
\item $F$ is full, and 
\item if a morphism $h$ in $\cA$ factors as $F(h)= fg$ in $\cA'$, with either $f$ or $g$ lifting to $\tilde{f}$ or $\tilde{g}$ in $\cA$, then the other lifts uniquely in such a way that $\tilde{f}\tilde{g}=h$.
\end{enumerate}

Since $\pi_n\C \to \C_0$ has a canonical section (given by $0 \in \pi_n$), it is full. 
Now, as $\pi_0\C$ is a groupoid, for all $f \in \C_0(x,y)$ the map $f_*:\pi_n\EEnd_{\C}(x) \to \pi_n(\HHom_{\C}(x,y), f)$ is an isomorphism. It therefore gives rise to an isomorphism  $\ad_f: \pi_n\EEnd_{\C}(x)\to \pi_n\EEnd_{\C}(y)$. We can write $\tilde{f}= f_*a$, $\tilde{g}=g_*b$ and $h= f_*g_*c$, so $\tilde{f}$ and $\tilde{g}$ determine each other by the formula $a= \ad_g(c b^{-1})$. 

Now, (F2) implies that $B\pi_0\C \to B\pi_0\cD$ is a Kan fibration, while (F1) implies that $B\C \to B\cD$ is a pointwise fibration. Therefore the Bousfield--Friedlander Theorem (\cite{sht} Theorem IV.4.9) gives that
\[
\diag B (\C\by_{\cD}\cB) \to  (\diag B\C)\by^h_{(\diag B\cD)}(\diag B\cB)
\]
is a weak equivalence. Finally, \cite{CRdiag}  shows that the canonical natural transformation
$\diag B\cA \to \bar{W}\cA$ is a weak equivalence for all $\cA$, which completes the proof.
\end{proof}

\begin{definition}
Given a category $\C$, define $c(\C)$ to be the set of isomorphism classes in $\C$.
\end{definition}

\begin{corollary}\label{piscat}
If $\C \in s\Cat$ is such that $\pi_0\C$ is a groupoid, then the homotopy fibre of $\bar{W}\C \to c(\pi_0\C)$ 
over $x \in \Ob \C$ is just
\[
\bar{W}\EEnd_{\C}(x).
\]
Thus $\pi_n(\bar{W}\C,x)= \pi_{n-1}\EEnd_{\C}(x)$, where the basepoint of the $\EEnd_{\C}(x)$ is taken to be the identity.
\end{corollary}

\begin{proof}
By using the model structure of \cite{bergner} Theorem 1.1, there is a morphism $g:\C \to \C'$ over $c(\C)$, with $c(\pi_0\C)= c(\pi_0\C')$, $\pi_*\HHom_{\C'}(gx,gy)\cong \pi_*\HHom_{\C}(x,y)$, and $\C' \to c(\C)$ a fibration. Replacing $\C$ with $\C'$, we may therefore assume that $\HHom_{\C}(x,y) $ is fibrant for all $x,y$.

Let $\cD$ be the groupoid with objects $\Ob \C$, and a unique morphism $x \to y$ for any pair of objects with the same image in $c(\pi_0\C)$. There is therefore a canonical morphism $\C \to \cD$, and this is a fibration. Therefore the proposition gives
\[
\bar{W}(\C\by_{\cD}\{x\}) \simeq (\bar{W} \C)\by^h_{(\bar{W}\cD)}\{x\}.
\]
Since $\bar{W}\cD \simeq c(\pi_0\C)$, the right-hand side is just the homotopy fibre of $\bar{W}\C \to c(\pi_0\C)$. Meanwhile, $\C\by_{\cD}\{x\} $ is just the simplicial category $\EEnd_{\C}(x)$ with unique object $x$ and morphisms $\HHom_{\C}(x,x)$, as required.
\end{proof}

\begin{corollary}\label{wscat}
If $\C$ and $\cD$ are simplicial categories for which $\pi_0\C$ and $\pi_0\cD$ are groupoids, then   a morphism $f:\C \to \cD$ is an $\infty$-equivalence if and only if
\[
\bar{W} \C \to \bar{W}\cD
\] 
is a weak equivalence.
\end{corollary}
\begin{proof}
We have just seen  that $f$ is  a $\bar{W}$-equivalence if and only if $c(\pi_0\C) \to c(\pi_0\cD)$ is an isomorphism, and $\pi_{n}\EEnd_{\C}(x) \to \pi_{n}\EEnd_{\cD}(fx) $ is an isomorphism for all $n \ge 0$. This is precisely the same as saying that $f$ is an $\infty$-equivalence.
\end{proof}

\begin{corollary}\label{hbycalc}
If $ \C\to \cD$ and $\cB \to \cD$ are arbitrary morphisms in $s\Cat$, with $\pi_0\C$ and $\pi_0\cD$ groupoids, then
\[
\bar{W}(\C\by_{\cD}^h\cB)\to (\bar{W} \C)\by^h_{(\bar{W}\cD)}(\bar{W}\cB) 
\]
is a weak equivalence, where the homotopy fibre product $\C\by_{\cD}^h\cB$ is taken in the model structure of \cite{bergner}.

In particular, if $\pi_0\cB$ is also a groupoid, then the fundamental groupoid $\pi_f ((\bar{W} \C)\by^h_{(\bar{W}\cD)}(\bar{W}\cB))$ is equivalent to the $2$-fibre product
\[
(\pi_0\C)\by^{(2)}_{(\pi_0\cD)}(\pi_0\cB),
\]
while  for any object $\alpha$ of $\C\by^{(2)}_{\cD}\cB $  given by $x \in \Ob \C, y \in \Ob \cB$ and $f \in \cD_0(\bar{x},\bar{y})$, we have
\[
\pi_n((\bar{W} \C)\by^h_{(\bar{W}\cD)}(\bar{W}\cB), \alpha) \cong \pi_{n-1}(\EEnd_{\C}(x)\by^h_{f_*, \HHom_{\cD}(\bar{x}, \bar{y}), f^* }\EEnd_{\cB}(y)).  
\]
\end{corollary}
\begin{proof}
Factorise $\C \to \cD$ as  $\C \to \C'\to \cD$, a weak equivalence followed by a fibration. By Corollary \ref{wscat}, 
\[
(\bar{W} \C)\by^h_{(\bar{W}\cD)}(\bar{W}\cB) \to (\bar{W} \C')\by^h_{(\bar{W}\cD)}(\bar{W}\cB)
\] 
is a weak equivalence. By Proposition \ref{fibscat}, the latter is weakly equivalent to $\bar{W}(\C'\by_{\cD}\cB)$, which is just $\bar{W}(\C\by^h_{\cD}\cB)$, since $s\Cat$ is right proper.

The description of the fundamental groupoid follows from the observation that for $\cA \in s\Cat$ and a discrete groupoid $\Gamma$, there are canonical isomorphisms 
\[
\Hom_{\bS}(\bar{W}\cA, B\Gamma) \cong \Hom_{s\Cat}(\cA, \Gamma) \cong \Hom_{\Cat}(\pi_0\cA, \Gamma). 
\]
For $\C'\to \cD$ a fibration, the map 
\[
\pi_0( \C'\by_{\cD}\cB)= (\pi_0\C')\by_{(\pi_0\cD)}(\pi_0\cB) \to   (\pi_0\C')\by^{(2)}_{(\pi_0\cD)}(\pi_0\cB)
\] 
is an equivalence, and the latter is equivalent to $(\pi_0\C)\by^{(2)}_{(\pi_0\cD)}(\pi_0\cB) $.  

Finally, for $\alpha:= (x,f,y)$ as above, also write $x$ for the image of $x$ in $\C'$. Since $\C' \to \cD$ lifts homotopy equivalences, there is a morphism $g: x \to z$ in $\C'$ mapping to $f$ in $\cD$. Then $\alpha$ defines the same element in $c(\pi_0\C')\by^{(2)}_{(\pi_0\cD)}(\pi_0\cB) $ as $\alpha':=(z, \id, y)$, so we need only describe homotopy groups based at $\alpha' $.

By Corollary \ref{piscat}, 
\[
\pi_n(\bar{W}(\C'\by_{\cD}\cB), \alpha')= \pi_{n-1}( \EEnd_{\C'}(z)\by_{\EEnd_{\cD}(\bar{y})}\EEnd_{\cB}(y), \alpha'),
\]
and we note that $\EEnd_{\C'}(z) \to \EEnd_{\cD}(\bar{y})$ is a Kan fibration, since $\C' \to \cD$ is a fibration. Now, since $g$ and $f$ are homotopy equivalences, they induce weak equivalences 
\begin{eqnarray*}
g^*: \EEnd_{\C'}(z) &\to& \HHom_{\C'}(x,z)\\
g_*:\EEnd_{\C'}(x)&\to& \HHom_{\C'}(x,z)\\
f^*:\EEnd_{\cD}(\bar{y})&\to& \HHom_{\cD}(\bar{x},\bar{y}),
\end{eqnarray*}
 giving  weak equivalences
\[
 \EEnd_{\C'}(z)\by_{\EEnd_{\cD}(\bar{y})}\EEnd_{\cB}(y) \to \HHom_{\C'}(x,z)\by^h_{\HHom_{\cD}(\bar{x},\bar{y}) }\EEnd_{\cB}(y) \la  \EEnd_{\C'}(x)\by^h_{f_*,\HHom_{\cD}(\bar{x},\bar{y}) }\EEnd_{\cB}(y).
\]  
Finally, the weak equivalence $\C \to \C'$ gives a weak equivalence $\EEnd_{\C}(x)\to\EEnd_{\C'}(x)$, yielding the expression required. 
\end{proof}

\section{Derived moduli from simplicial categories}\label{modsimpcat}

\subsection{Background on representability}

Denote the  category of simplicial rings by $s\Ring$, the category of simplicial $R$-algebras by $s\Alg_R$, and the category of simplicial $R$-modules by $s\Mod_R$. If $R$ is a (dg) $\Q$-algebra, we let $dg_+\Alg_R$ be the category of differential graded-commutative $R$-algebras in non-negative chain degrees, and $ dg_+\Mod_R$ the category of $R$-modules in  chain complexes in non-negative chain degrees.  

\begin{definition}\label{normdef}
 Given a simplicial abelian group $A_{\bt}$, we denote the associated normalised chain complex  by $NA$. Recall that this is given by  $N(A)_n:=\bigcap_{i>0}\ker (\pd_i: A_n \to A_{n-1})$, with differential $\pd_0$. Then $\H_*(NA)\cong \pi_*(A)$.

Using the Eilenberg--Zilber shuffle product,  normalisation  $N$ extends to a functor
$$
N:s\Alg_R \to dg_+\Alg_R.
$$ 
When $R$ is a $\Q$-algebra, this functor is a right Quillen equivalence.
\end{definition}

\begin{definition}
Define $dg_+\cN_R$ (resp. $s\cN_R$)  to be the full subcategory of $dg_+\Alg_R$ (resp. $s\Alg_R$) consisting of objects $A$
 for which the map $A \to \H_0A$ (resp. $A \to \pi_0A$) has nilpotent kernel. Define $dg_+\cN_R^{\flat}$ (resp. $s\cN_R^{\flat}$) to be the full subcategory of $dg_+\cN_R$ (resp. $s\cN_R$) consisting of objects $A$
 for which $A_i=0$ (resp. $N_iA=0$) for all $i \gg 0$.
\end{definition}

From now on, we will write $d\cN^{\flat}$ (resp. $d\Alg_R$, resp. $d\Mod_R$) to mean either $s\cN_R^{\flat}$ (resp. $s\Alg_R$, resp. $s\Mod_R$) or $dg_+\cN_R^{\flat}$ (resp. $dg_+\Alg_R$, resp. $dg_+\Mod_R$), noting that we only use chain algebras   in characteristic $0$. For $A \in d\Alg_R$, write $I_A:= \ker(A \to \H_0A)$.

\begin{definition}
Say that a  surjection $A \to B$ in  $dg_+\Alg_R$ (resp. $s\Alg_R$) is a \emph{tiny acyclic extension} if the kernel $K$ satisfies $I_A\cdot K=0$, and  $K$ (resp. $NK$) is of the form $\cone(M)[-r]$ for some $\H_0A$-module  (resp. $\pi_0A$-module) $M$. In particular, $\H_*K=0$. 
\end{definition}

\begin{definition}
Say that a functor $F:  d\cN^{\flat} \to \C$ to a model category $\C$ is homotopy-preserving if for all tiny acyclic extensions $A \to B$, the map
\[
F(A) \to F(B)
\]
is a weak equivalence.
\end{definition}

\begin{definition}\label{hhgsdef}
We say that a 
functor 
$$
F:d\cN^{\flat} \to \C
$$
to a model category $\C$ is homotopy-homogeneous if  for all square-zero extensions $A \to B$ and all maps $C \to B$ in $d\cN^{\flat} $, the natural map
$$
F(A\by_BC) \to F(A)\by^h_{F(B)}F(C)
$$
to the homotopy fibre product is a weak  equivalence.
\end{definition}

\begin{definition}\label{Tdef}
Given a homotopy-preserving homotopy-homogeneous  functor $F: d\cN^{\flat} \to \bS$,  $A \in d\cN^{\flat}$, and a point $x \in F(A)$,  define
$$
T_x(F/R): d\Mod_A \to \bS
$$
by
$$
T_x(F/R)(M):= F(A\oplus M)\by^h_{F(A)}\{x\}.
$$
\end{definition}

\begin{definition}\label{totcohodef} 
Given $F,A,x$ as above, and $M \in d\Mod_A$,  define
$$
\DD^{n-i}_x(F/R,M):= \pi_i (T_x(F/R) (M[-n])),
$$
noting that this is well-defined  by \cite{drep} Lemma \ref{drep-adf}.
\end{definition}

\begin{proposition}\label{obs}
If $F:s\Alg_R \to \bS$ is homotopy-preserving and homotopy-homogeneous, then   
  for any square-zero extension $e:I \to A \xra{f} B$ in $\C$, there is a sequence of sets
$$
\pi_0(FA)\xra{f_*} \pi_0(FB) \xra{o_e} \Gamma(FB,\DD^1(F/R, I)), 
$$  
where $\Gamma(-)$ denotes the global section functor. 
This is exact in the sense that the fibre of $o_e$ over $0$ is the image of $f_*$. 
 Moreover,  there is a group action of $\DD^0_x(F/R, I)$ on the fibre of $\pi_0(FA) \to \pi_0(FB)$ over $x$, whose orbits are precisely the fibres of $f_*$. 

For any $y \in F_0A$, with $x=f_*y$, the fibre of $FA \to FB$ over $x$ is weakly equivalent to $T_{x}(F/R,I)$, and the sequence above 
extends to a long exact sequence
$$\xymatrix@R=0ex{
\cdots \ar[r]^-{e_*} &\pi_n(FA,y) \ar[r]^-{f_*}&\pi_n(FB,x) \ar[r]^-{o_e}& \DD^{1-n}_{y}(F/R,I) \ar[r]^-{e_*} &\pi_{n-1}(FA,y)\ar[r]^-{f_*}&\cdots\\ &\cdots \ar[r]^-{f_*}&\pi_1(FB,x) \ar[r]^-{o_e}& \DD^0_{y}(F/R,I)  \ar[r]^-{-*y} &\pi_0(FA).
}
$$
\end{proposition}
\begin{proof}
 This is \cite{drep} Proposition \ref{drep-obs}
\end{proof}

\begin{definition}
 Given a functor  $F:  d\cN_R^{\flat} \to \bS$, define $\pi^0F: \Alg_{\pi_0R} \to \bS$ by $\pi^0F(A)=F(A)$.
\end{definition}

The following variant of Lurie's Representability Theorem is \cite{drep} Theorem \ref{drep-lurierep3}:
\begin{theorem}\label{lurierep3}
Let $R$ be a derived G-ring admitting a dualising module (in the sense of \cite{lurie} Definition 3.6.1) and   take a functor $F: d\cN_R^{\flat} \to \bS$. Then $F$ is the restriction of  an almost finitely presented geometric derived $n$-stack $F':d\Alg_R \to \bS$  if and only if 
 the following conditions hold

\begin{enumerate}
 
\item $F$ is homotopy-preserving.

\item For all discrete rings $A$, $F(A)$ is $n$-truncated, i.e. $\pi_iF(A)=0$ for all $i>n$ .

\item\label{cohesive} 
$F$ is homotopy-homogeneous.

\item $\pi^0F:\Alg_{\pi_0R} \to \bS$ is a hypersheaf for the \'etale topology. 

\item\label{afp1a} $\pi_0\pi^0F:  \Alg_{\pi_0R} \to \Set$  preserves filtered colimits.

\item\label{afp1b} For all $A \in \Alg_{\pi_0R}$ and all $x \in F(A)$, the functors $\pi_i(\pi^0F,x): \Alg_A \to \Set$  preserve filtered colimits for all $i>0$.

\item 
for all finitely generated integral domains $A \in \Alg_{\pi_0R}$, all $x \in F(A)_0$ and all \'etale morphisms $f:A \to A'$, the maps
\[
\DD_x^*(F, A)\ten_AA' \to \DD_{fx}^*(F, A')
\]
are isomorphisms.

\item\label{afp2} for all finitely generated $A \in \Alg_{\pi_0R}$  and all $x \in F(A)_0$, the functors $\DD^i_x((F/R), -): \Mod_A \to \Ab$ preserve filtered colimits for all $i>0$.

\item for all finitely generated integral domains $A \in \Alg_{\pi_0R}$  and all $x \in F(A)_0$, the groups $\DD^i_x(F/R, A)$ are all  finitely generated $A$-modules.

\item for all complete discrete local Noetherian  $\pi_0R$-algebras $A$, with maximal ideal $\m$, the map
$$
\pi^0F(A) \to {\Lim}^h F(A/\m^r)
$$
is a weak equivalence.
\end{enumerate}
Moreover, $F'$ is uniquely determined by $F$ (up to weak equivalence).
\end{theorem}

\begin{remark}\label{formalexistrk}
The Milnor exact sequence (\cite{sht} Proposition 2.15) gives a  sequence
\[
\bt \to \Lim^1_r \pi_{i+1}F(A/\m^r) \to \pi_i({\Lim}^h F(A/\m^r))  \to \Lim_r \pi_iF(A/\m^r)\to \bt,
\]
which is exact as groups for $i \ge 1$ and as pointed sets for $i=0$. Thus the final condition above can be rephrased to say that the map 
\[
f_0:\pi_0F(A) \to \Lim_r  \pi_0F(A/\m^r)
\] 
is surjective,  that for all $x \in F(A)$ the maps 
\[
f_{i,x}:\pi_i(FA,x) \to \Lim_r  \pi_i(F(A/\m^r),x)
\] 
are surjective for all $i \ge 1$ and that the resulting maps
\[
\ker f_{i,x} \to  \Lim^1_r \pi_{i+1}(F(A/\m^r),x)
\]
are surjective for all $i \ge 0$.

Now, we can say that an inverse system $\{G_r\}_{r \in \N}$ of groups satisfies the Mittag-Leffler condition if for all $r$, the images $\im(G_s \to G_r)_{s\ge r}$ satisfy the descending chain condition. In that case, the usual abelian proof (see e.g. \cite{W} Proposition 3.5.7) adapts to show that $\Lim^1 \{G_r\}_r = 1$. 

Hence, if each system $\{\im( \pi_{i}(F(A/\m^s),x)\to  \pi_{i}(F(A/\m^r),x))\}_{s \ge r}$ satisfies the Mittag-Leffler condition for $i \ge 1$, then the final condition of the theorem
 reduces to requiring that the maps
\[
\pi_iF(A) \to \Lim_r  \pi_iF(A/\m^r)
\]
be isomorphisms for all $i$.
\end{remark}

\begin{remark}\label{cflurie}
Note that there are slight differences in terminology between \cite{hag2} and \cite{lurie}. In the former, only disjoint unions of affine schemes are $0$-representable, so arbitrary schemes are $2$-geometric stacks, and Artin stacks are $1$-geometric stacks if and only if they have affine diagonal. In the latter, algebraic spaces are $0$-stacks.  A geometric $n$-stack  is called $n$-truncated in \cite{hag2}, and it follows easily that every $n$-geometric stack in \cite{hag2} is $n$-truncated. Conversely, every  geometric $n$-stack  is  $(n+2)$-geometric. 

Theorem \ref{lurierep3} takes the convention from \cite{lurie}, so ``geometric derived $n$-stack'' means ``$n$-truncated derived geometric stack''.
\end{remark}

\subsection{Formal quasi-smoothness}

\begin{definition}
Say that a natural transformation $\eta: F \to G$ of functors $F, G:  d\cN^{\flat} \to \C$ to a model category $\C$ is homotopic if for all tiny acyclic extensions $A \to B$, the map
\[
F(A) \to F(B)\by_{G(B)}G(A)
\]
is a trivial fibration. Say that $F$ is homotopic if $F \to \bt$ is so, where $\bt$ denotes the final object of $\C$.
\end{definition}

\begin{definition}
Say that a natural transformation $\eta: F \to G$ of functors $F, G:  d\cN^{\flat} \to \C$ to a model category $\C$ is formally quasi-presmooth  if for all square-zero extensions $A \to B$, the map
\[
F(A) \to F(B)\by_{G(B)}G(A)
\]
is a fibration.

Say that $\eta$ is formally quasi-smooth if it is formally quasi-presmooth  and homotopic. 
\end{definition}

\begin{definition}
Say that a natural transformation $\eta: F \to G$ of functors $F, G:  d\cN^{\flat} \to \bS$ is formally presmooth (resp. formally smooth) if it is formally quasi-presmooth (resp. formally quasi-smooth) and
for all square-zero extensions $A \to B$, the map
\[
F(A) \to F(B)\by_{G(B)}G(A)
\]
is surjective on $\pi_0$.
\end{definition}

\begin{definition}
Say that a natural transformation $\eta: F \to G$ of functors on $  d\cN^{\flat}$ is formally \'etale if 
for all square-zero extensions $A \to B$, the map
\[
F(A) \to F(B)\by_{G(B)}G(A)
\]
is an isomorphism.
\end{definition}

\subsection{Homogeneity}

\begin{definition}
Say that a natural transformation $F \to G$ of functors on $ d\cN^{\flat}$ is (relatively) homogeneous if for all square-zero extensions $A \to B$, the map
\[
F(A\by_BC ) \to G(A\by_BC)\by_{[G(A)\by_{G(B)}G(C)]}[F(A)\by_{F(B)}F(C)]
\]
is an isomorphism. Say that $F$ is homogeneous if $F \to \bt$ is relatively homogeneous.
\end{definition}

\begin{proposition}\label{fethgs}
Let $\alpha:F \to G$ be a formally \'etale morphism of functors $F, G :  d\cN^{\flat} \to \Set$. If $G$ is homogeneous, then so is $F$. Conversely, if  $\alpha$ is surjective  and $F$ is  homogeneous, then so is $G$.
\end{proposition}
\begin{proof}
Take a square-zero extension $A \to B$, and a morphism $C \to B$, noting that $A\by_BC \to C$ is then another square-zero extension. 
Since $\alpha$ is  formally \'etale,
\[
F(A\by_BC)\cong G(A\by_BC)\by_{GC}FC, 
\]
and
\begin{eqnarray*}
FA\by_{FB}FC&\cong&  [GA\by_{GB}FB]\by_{FB}FC\\
&=&  GA\by_{GB}FC \\
&=& [GA\by_{GB}GC]\by_{GC}FC.
\end{eqnarray*}
Thus homogeneity of $G$ implies homogeneity of $F$, and if $\pi_0FC \to \pi_0GC$ is surjective for all $C$, then homogeneity of $F$ implies homogeneity of $G$.
\end{proof}

\subsection{Simplicial categories}

\begin{definition}\label{2fibre}
Given functors $\cA \xra{F} \cB \xla{G} \C$ between categories, define the 2-fibre product $\cA\by^{(2)}_{\cB}\C$ as follows. Objects of $\cA\by^{(2)}_{\cB}\C  $ are triples $(a,\theta, c)$, for $a \in \cA, c \in \C$ and $\theta: Fa \to Gc$ an isomorphism in $\cB$. A morphism in   $\cA\by^{(2)}_{\cB}\C$ from $(a,\theta, c)$ to $(a',\theta', c')$ is a pair $(f,g)$, where $f:a \to a'$ is a morphism in $\cA$ and $g: c\to c'$ a morphism in $\C$, satisfying the condition that
\[
Gg \circ \theta = \theta' \circ Ff.
\] 
\end{definition}

\begin{remark}
This definition has the property  that $\cA\by^{(2)}_{\cB}\C$ is a model for the 2-fibre product in the 2-category of categories. However, we will always use the notation $\cA\by^{(2)}_{\cB}\C$ to mean the specific model of Definition \ref{2fibre}, and not merely any equivalent category. 

Also note that
\[
\cA\by^{(2)}_{\cB}\C= (\cA\by^{(2)}_{\cB}\cB)\by_{\cB}\C
\]
\end{remark}

\begin{definition}\label{s2fibre}
Given functors $\cA \xra{F} \cB \xla{G} \C$ between simplicial categories, define the 2-fibre product $\cA\by^{(2)}_{\cB}\C$ as follows. Objects of $\cA\by^{(2)}_{\cB}\C  $ are triples $(a,\theta, c)$, for $a \in \cA, c \in \C$ and $\theta: Fa \to Gc$ an isomorphism in $\cB_0$. Morphisms are given by setting
\[
\HHom_{\cA\by^{(2)}_{\cB}\C}( (a,\theta, c), (a',\theta', c')):=  
\{(f,g) \in  \HHom_{\cA}(a,a') \by \HHom_{\C}(c,c') \,:\,  Gg \circ \theta = \theta' \circ Ff \}.
\]
\end{definition}

\begin{definition}\label{2fibrn}
Say that a morphism $F:\cA \to \cB$ in $s\Cat$ is a (trivial) 2-fibration if $\cA\by_{\cB}^{(2)}\cB \to \cB$ is a (trivial) fibration in the sense of Definition \ref{scatfibdef}. Explicitly, $F$ is a $2$-fibration if
\begin{enumerate} 
\item[(F1)] for any objects $a_1$ and $a_2$ in $\cA$, the map
$\HHom_{\cA}(a_1, a_2)\to \HHom_{\cB}(Fa_1, Fa_2)$ 
is a fibration of simplicial sets;
\item[(F2)] for any objects $a_1 \in\cA$, $b \in \cB$, and any homotopy equivalence $e :
Fa_1 \to b$ in $\cB$, there is an object $a_2 \in \cA$,  a homotopy equivalence
$d : a_1 \to a_2$ in $\cA$ and an isomorphism $\theta: Fa_2 \to b$   such that $\theta\circ Fd = e$.
\end{enumerate}
Likewise, $F$ is a trivial $2$-fibration if
\begin{enumerate} 
\item[(TF1)] for any objects $a_1$ and $a_2$ in $\cA$, the map
$\HHom_{\cA}(a_1, a_2)\to \HHom_{\cB}(Fa_1, Fa_2)$ 
is a trivial fibration of simplicial sets;
\item[(TF2)] $F_0: \cA_0 \to \cB_0$ is essentially surjective on objects.
\end{enumerate}
\end{definition}

\begin{definition}
Say that a natural transformation $\eta: \cD \to \cE$ of functors $\cD:  d\cN^{\flat} \to s\Cat$ is 2-homotopic  if for all tiny acyclic extensions $A \to B$, the map
\[
\cD(A) \to \cD(B)\by^{(2)}_{\cE(B)}\cE(A)
\]
is a trivial 2-fibration. Say that $\cD$ is 2-homotopic if $\cD \to \bt$ is so.
\end{definition}

\begin{lemma}\label{2htpicgood}
If $\cD$ is 2-homotopic, then $\cD$ is homotopy-preserving.
\end{lemma}
\begin{proof}
By \cite{drep} Lemma \ref{drep-wtiny} and the proof of \cite{drep} Proposition \ref{drep-cNhat}, it suffices to show that $\cD$ maps tiny acyclic extensions $A \to B$  to weak equivalences. Since $\cD_0(A)\to \cD_0(B)$ is essentially surjective on objects, so is $\pi_0\cD(A) \to \pi_0\cD(B)$. On morphisms, we have $\HHom_{\cD(A)}(x,y) \to \HHom_{\cD(B)}(x,y)$ a trivial fibration (and hence a weak equivalence) for all $x, y \in \cD(A)$, so $\cD(A) \to \cD(B)$ is indeed a weak equivalence.
\end{proof}

\begin{definition}
Say that a natural transformation $\eta: \cD \to \cE$ of functors $\cD,\cE:  d\cN^{\flat} \to s\Cat$ is formally 2-quasi-presmooth  if for all square-zero extensions $A \to B$, the map
\[
\cD(A) \to \cD(B)\by^{(2)}_{\cE(B)}\cE(A)
\]
is a 2-fibration. If $\eta$ is also 2-homotopic, then we say that $\eta$ is formally 2-quasi-smooth.

Say that $\cD$ is formally 2-quasi-presmooth (resp. formally 2-quasi-smooth) if $\cD \to \bt$ is so.
\end{definition}

\begin{definition}
Say that a functor $\cF :  d\cN^{\flat} \to s\Cat$ is  2-homogeneous if for all square-zero extensions $A \to B$, the map
\[
\cF(A\by_BC ) \to \cF(A)\by_{\cF(B)}^{(2)}\cF(C)
\]
is essentially surjective on objects (for $\cF_0$), and an isomorphism on $\HHom$-spaces.
\end{definition}

\begin{lemma}\label{2hgsgood}
If $\cD:  d\cN^{\flat} \to s\Cat$ is 2-homogeneous and formally 2-quasi-presmooth, 
then
$\cD$ is homotopy-homogeneous.
\end{lemma}
\begin{proof}
Take  a square-zero extension $A \to B$ and a map $C \to B$. We first need to show that
\[
\theta: \pi_0\cD(A\by_BC) \to \pi_0(\cD(A)\by^h_{\cD(B)}\cD(C))
\]
is essentially surjective on objects. 

Any object of the right-hand side is represented by   objects $x \in \cD(A)$, $y \in \cD(B)$, and a homotopy equivalence $f:\bar{x}\to \bar{y}$ in $\cD(B)$. Since $\cD$ is formally 2-quasi-presmooth, we may lift $f$ to some map $\tilde{f}: x \to z$ in $\cD(A)_0$ with $g:\bar{z} \cong \bar{y}$ in $\cD(B)_0$. Since $\cD$ is 2-homogeneous, there must then be an object $t \in \cD(A\by_BC) $ whose image in $\cD(A)_0\by_{\cD(B)_0}^{(2)}\cD(C)_0$ is isomorphic to $(z,g,y)$. Since $\tilde{f}$ gives a homotopy equivalence from this to $(x, f, y)$, we have shown that $\theta$ is essentially surjective on objects.
 
For morphisms, we have
\begin{eqnarray*}
\HHom_{\cD(A\by_BC)}(t,t')&\cong& \HHom_{\cD(A)}(t,t')\by_{\HHom_{\cD(B)}(t,t')}\HHom_{\cD(C)}(t,t')\\
&\simeq& \HHom_{\cD(A)}(t,t')\by^h_{\HHom_{\cD(B)}(t,t')}\HHom_{\cD(C)}(t,t'),
\end{eqnarray*}
the first isomorphism following from homogeneity, and the second because $\HHom_{\cD(A)}(t,t')\to {\HHom_{\cD(B)}(t,t')} $ is a fibration, by quasi-smoothness. This completes the proof that
\[
\cD(A\by_BC) \simeq \cD(A)\by^h_{\cD(B)}\cD(C)
\]
in $s\Cat$.
\end{proof}

\begin{definition}
Say that $\cD:  d\cN^{\flat} \to s\Cat$ has formally quasi-presmooth (resp. formally quasi-smooth) $\HHom$-spaces if $\HHom_{\cD}(x,y): A \da   d\cN^{\flat} \to \bS$ is formally quasi-presmooth (resp. formally quasi-smooth) for all  $A \in  d\cN^{\flat}$ and $x,y \in \cD(A)$.
\end{definition}

\begin{definition}
Given a simplicial category $\C$, define  $\cW(\C)$ to be the full  simplicial subcategory in which morphisms are maps whose images in $\pi_0\C$ are invertible. In particular, this means that $\pi_0\cW(\C)$ is the core (maximal subgroupoid) of $\pi_0\C$.
\end{definition}

\begin{proposition}\label{wetale}
If $\cD:  d\cN^{\flat} \to s\Cat$ is 2-homogeneous with formally quasi-smooth $\HHom$-spaces,
then $\cW(\cD) \to \cD$ is formally \'etale, in the sense that
for any square-zero extension $A \to B$, the map
\[
\cW(\cD)(A) \to \cD(A)\by_{\cD(B)}\cW(\cD)(B)
\]
is an isomorphism.
\end{proposition}
\begin{proof}
Since $f \in \HHom_{\C}(x,y)_n$ lies in $\cW(\C)$ if and only if $(\pd_0)^nf$ is a homotopy equivalence, 
it suffices to show that
\[
\cW(\cD)_0(A) \to \cD_0(A)\by_{\cD_0(B)}\cW(\cD)_0(B)
\]
is an isomorphism.

Take a morphism   $f \in \HHom_{\cD(B)}(x,y)_0$, with objects $x,y$ lifting to $x,y\in\cD(A)$.  If we set $\DD^n_f(M):= \DD^{n}_f(\HHom_{\cD}(x,y)/A,M) $, then  Proposition \ref{obs} shows that the obstruction $o(f)$ to lifting $\bar{f}$ to  $\HHom_{\cD(A)}(x,y)_0$ lies in $\DD^1_{\bar{f}}(I)$, where $I= \ker(A \to B)$. 

Assume that $f$  is a  homotopy equivalence, so there exist $g \in \HHom_{\cD(B)}(y,x)_0$ and homotopies $h^x \in \HHom_{\cD(B)}(x,x)_1, h^y \in \HHom_{\cD(B)}(y,y)_1 $ such that 
\[
        \pd_1h^x= \id_x,\quad \pd_1 h^y=\id_y, \quad\pd_0h^x= g\circ \bar{f},\quad \pd_0h^y= \bar{f}\circ g.
\]
 We now need to study the construction of obstruction maps.

For the acyclic extension $\tilde{B} \to B$ as in the Proof of Proposition \ref{obs}, we may lift $f$ to $\tilde{f}$ in $\cD(\tilde{B})$, and do the same for $g$. Then the obstruction $o(f\circ g)$ to lifting $f\circ g$ is the class of the image of $\tilde{f}\circ \tilde{g}$ in $\pi_0\HHom{\cD(B \oplus I[-1])}(y,y)$, which is just
\[
o(f\circ g)= f\circ g+ f_*o(g) +g^*o(f) \in \DD_{f \circ g}^1(I).
\]
Since $\HHom_{\cD(\tilde{B})}(x,x)\to \HHom_{\cD(B)}(x,x)$ is a trivial fibration (by quasi-smoothness), we may also lift the homotopies $h^x, h^y$ to give  homotopies $\tilde{h}^x\in\HHom_{\cD(\tilde{B})}(x,x)_1, \tilde{h}^y\in \HHom_{\cD(\tilde{B})}(y,y)_1$ making $\tilde{f}$ and $\tilde{g}$ homotopy inverses. Then  $\tilde{h}^y$ ensures that $o(f \circ g)=0 $ (since the image of $\tilde{f}\circ \tilde{g}$ is homotopic to $\id_y$), and therefore $f_*o(g)= -g^*o(f) \in \DD^1_y(I)$. Likewise $g_*o(f)= -f^*o(g) \in \DD^1_x(I)$. Since $o(g)= (gf)_*(fg)^*o(g)$, this gives $o(g)= -g_*g^*o(f)$

If $f$ lifts to $f'\in\HHom_{\cD(A)}(x,y)_0$, then $o(f)=0$, and the formula above ensures that $o(g)=0$, so $g$ also lifts to $g'$ in  $\cD(A)$. Now, $f' \circ g'$ is in the fibre of $\EEnd_{\cD(A)}(y)_0 \to \EEnd_{\cD(B)}(y)_0$ over $\id_y$, so $[f' \circ g'] \in \pi_0\EEnd_{\cD(A)}(y)$ lies in the image of 
\[
\pi_0T_{\id_y}( \EEnd_{\cD}(y)/A, I)  = \DD^0_y(I). 
\]
Another choice for $g'$ is of the form $g' + a$ for $a \in \DD^0_g(I)$, and then $[f'\circ (g'+1)] = [f' \circ g']+ f_*a$. Since $f$ is a homotopy equivalence, $f_*: \DD^0_g(I)\to \DD^0_y(I)$ is an isomorphism, so we may set $a= -(f_*)^{-1}[f'\circ g']$. Without loss of generality, we may therefore assume that $f'$ and $g'$ are homotopy inverses. 

Finally, write $[g' \circ f']= \id_x +b $, for $b \in \DD^0_x(I)$. Since $f'\circ g'$ is homotopic to $\id_y$, applying $f_*$ to $b$ gives 
\[
[f']+ f_*b = [f'\circ g' \circ f'] = [f'], 
\]
so $f_*b=0 \in \DD^0_{f}(I)$. Thus $b=0$, as $f_*: \DD^0_x(I)\to  \DD^0_{f}(I)$ is an isomorphism.
\end{proof}

\begin{corollary}\label{w2good}
If $\cD:  d\cN^{\flat}\to s\Cat$ is  formally 2-quasi-smooth and 2-homogeneous, then so is $\cW(\cD)$.
\end{corollary}
\begin{proof}
For any square-zero extension $A \to B$ in $  d\cN^{\flat} $, we have
\[
\cW(\cD)(A)= \cD(A)\by_{\cD(B)}\cW(\cD)(B),
\]
by Proposition \ref{wetale}. Thus
\[
\cW(\cD)(A)\by^{(2)}_{\cW(\cD)(B)}\cW(\cD)(B)= \cD(A)\by^{(2)}_{\cD(B)}\cW(\cD)(B),
\]
so $\theta: \cW(\cD)(A)\by^{(2)}_{\cW(\cD)(B)}\cW(\cD)(B) \to  \cW(\cD)(B) $ is the pullback of $  \cD(A)\by^{(2)}_{\cD(B)}\cD(B) \to \cD(B) $ along $\cW(\cD)(B) \to \cD(B)$. By 2-quasi-smoothness of $\cD$, this means that $\theta$ is a fibration in $s\Cat$, and a trivial fibration whenever $A \to B$ is a tiny acyclic extension. Thus $\cW(\cD)$ is formally 2-quasi-smooth.

To see that $\cW(\cD)$ is 2-homogeneous, first observe that $\cD(A)\by^{(2)}_{\cD(B)}\cD(C)$ and $\cW(\cD)(A)\by^{(2)}_{\cW(\cD)(B)}\cW(\cD)(C)$ have the same objects. For a square-zero extension $A \to B$ and a map $C \to B$, 2-homogeneity of $\cD$ then gives that
\[
\cW(\cD)_0(A\by_BC) \to \cW(\cD)_0(A)\by^{(2)}_{\cW(\cD)_0(B)}\cW(\cD)_0(C)
\]
is essentially surjective on objects. Proposition \ref{wetale} combines with 2-homogeneity of $\cD$ and Lemma \ref{fethgs} to give that $\HHom_{\cW(\cD)}(x,y)$ is homogeneous, since $\HHom_{\cW(\cD)}(x,y)\to \HHom_{\cD}(x,y) $ is formally \'etale.
\end{proof}
 


\begin{proposition}\label{descentlemma}
 Take a ring $R$, a functor $\cF \co \Alg_{R}\to s\Cat$ for which $\bar{W}\cW(\cF) \co \Alg_{R}\to \bS$ is  an \'etale hypersheaf, and a subfunctor $\cE\subset \cF$ for which each $\cE(A)$ is a full simplicial subcategory of $\cF(A)$. 

Then $\bar{W}\cW(\cE)$ is an \'etale hypersheaf if and only if for any $A \in \Alg_R$ and any  \'etale cover $\{f_i \co A \to B_i\}_{i \in I}$, the map
\[
 c\pi_0\cE(A) \to c\pi_0\cF(A)\by_{(\prod_ic\pi_0\cF(B_i))}(\prod_i c\pi_0\cE(B_i))
\]
is surjective, where $c\C$ denotes the set of isomorphism classes in a category $\C$.
\end{proposition}
\begin{proof}
The ``only if'' implication is almost automatic, so we will just prove the ``if'' implication.
First, take a finite set $\{A_i\}_{i \in I}$ of $R$-algebras. Since $F:= \bar{W}\cW(\cF) $ is a hypersheaf, the map
\[
 F(\prod_i A_i) \to \prod_iF(A_i)
\]
is a weak equivalence. 

Setting $E:= \bar{W}\cW(\cE)$, the fact that $\cE \subset \cF$ is a full simplicial subcategory implies that the map $E \to F$ is injective on $\pi_0$ and an isomorphism on all other homotopy groups. Writing $A:= \prod_i A_i$, it follows that
\[
 \pi_n(E(A), x)\cong \prod_i\pi_n(E(A_i), x)  
\]
for all $n>0$ and all $x$, and that $ \pi_0E(A)\into \prod_i\pi_0E(A_i)$.

Now $\pi_0E= c\pi_0\cE$ and the set $\{A \to A_i\}_{i \in I}$ is an \'etale cover, so we also have
\[
 \pi_0E(A) \onto \pi_0F(A)\by_{(\prod_i\pi_0F(A_i))}(\prod_i \pi_0E(A_i)) =  \prod_i \pi_0E(A_i),
\]
  and hence
\[
 E(A) \simeq \prod_iE(A_i).
\]

Now take an affine \'etale  hypercover $B \to B^{\bt}$; we need to show that the map
\[
 \phi\co E(B) \to \ho\Lim_n E(B^n)
\]
is a weak equivalence. Because $E \simeq \pi_0E\by_{\pi_0F}F$, we have
\begin{eqnarray*}
 \ho\Lim_n E(B^n) &\simeq& (\ho\Lim_n \pi_0E(B^n))\by_{(\ho\Lim_n\pi_0F(B^n))}(\ho\Lim_n F(B^n))\\
&\simeq& (\Lim_n \pi_0E(B^n))\by_{(\Lim_n\pi_0F(B^n))}(\ho\Lim_n F(B^n))\\
&\simeq& (\Lim_n \pi_0E(B^n))\by_{(\Lim_n\pi_0F(B^n))}F(B).
\end{eqnarray*}

Since $\pi_0E \subset \pi_0F$ and $\Lim_n \pi_0E(B^n) \subset \pi_0E(B^0)$, we can rewrite this as
\[
 \pi_0E(B^0)\by_{\pi_0F(B^0)}F(B),
\]
so $\phi$ induces an isomorphism on all higher homotopy groups and is injective on $\pi_0$.
Because $B \to B^0$ is an \'etale cover, the map
\[
 \pi_0E(B) \to \pi_0F(B)\by_{\pi_0F(B^0)}\pi_0E(B^0)
\]
is surjective, so $\phi$ is also surjective on $\pi_0$, and hence a weak equivalence.
 \end{proof}

\section{Derived moduli of schemes}\label{modschsn}
\subsection{Derived moduli of affine schemes}

\begin{definition}\label{dgcof}
Given $A \in dg_+\Alg_R$,  say that $B \in dg_+\Alg_A $ is quasi-free if  the underlying map $A_* \to B_*$   of graded algebras is freely generated. Say that $B \in  dg_+\Alg_A $ is cofibrant if it is a retract of a quasi-free dg algebra.

This is equivalent to satisfying the condition that 
 for all acyclic surjections $C' \to C$ in $dg_+\Alg_A$, the map $\Hom_{dg_+\Alg_A }(B,C') \to \Hom_{dg_+\Alg_A }(B,C)$ is surjective.
\end{definition}

\begin{definition}\label{scof}
Given $A \in s\Alg_R$,  say that $B \in s\Alg_A$ is quasi-free if   there are sets $\Sigma_q \subset B_q$ with $B_q= A_q[\Sigma_q]$, such that $\bigcup_q \Sigma_q$ is closed under the degeneracy operations.
Say that $B \in  s\Alg_A $ is cofibrant if it is a retract of a quasi-free simplicial algebra.

This is equivalent to satisfying the condition that 
 for all acyclic surjections $C' \to C$ in $s\Alg_A$, the map $\Hom_{s\Alg_A }(B,C') \to \Hom_{s\Alg_A }(B,C)$ is surjective.
\end{definition}

\begin{definition}\label{horndef}
Let $\Delta^n \in \bS$ be the standard $n$-simplex, and $\pd\Delta^n\in \bS$ its boundary.  Given $0 \le k \le n$, define the $k$th horn $\L^n_k$ of $\Delta^n$ to be the simplicial set obtained from $\Delta^n$ by removing the interior and the $k$th face. See \cite{sht} \S I.1 for explicit descriptions. 
\end{definition}

\begin{definition}
There is a simplicial structure on $s\Alg_A$ given as follows. For  $B\in s\Alg_A$ and $K \in \bS$,  $B^K$ is defined by 
$$
(B^K)_n:= \Hom_{\bS}(K \by \Delta^n, B).
$$ 
Spaces  $\HHom_{s\Alg_A }(C,B) \in \bS$ of morphisms are then given by
$$
\HHom_{s\Alg_A }(C, B)_n:= \Hom_{s\Alg_A }(C, B^{\Delta^n}).
$$
\end{definition}

\begin{definition}
As in \cite{drep} Lemma \ref{drep-simplicialstr}, put a  simplicial structure on $dg_+\Alg_A$  as follows.
First set $\Omega_n$ to be the  cochain algebra 
$$
\Q[t_0, t_1, \ldots, t_n,dt_0, dt_1, \ldots, dt_n ]/(\sum t_i -1, \sum dt_i)
$$  
of rational differential forms on the $n$-simplex $\Delta^n$.
 These fit together to form a simplicial complex $\Omega_{\bt}$ of DG-algebras, and we define $B^{\Delta^n}$ as the good truncation $B^{\Delta^n}:= \tau_{\ge 0}(B \ten \Omega_n)$.
Spaces  $\HHom_{dg_+\Alg_A }(C,B) \in \bS$ of morphisms are then given by
$$
\HHom_{dg_+\Alg_A }(C, B)_n:= \Hom_{dg_+\Alg_A }(C, B^{\Delta^n}).
$$
\end{definition}

\begin{definition}
For  $A \in dg_+\Alg_R$ (resp. $A \in s\Alg_R$),
define $\ALG(A)$ to be the full simplicial subcategory of the category $dg_+\Alg_A $ (resp. $s\Alg_A$) on  cofibrant objects.
\end{definition}

\begin{proposition}\label{alghgs}
The functor $\ALG:  d\cN^{\flat} \to s\Cat$ defined above is 2-homogeneous and formally 2-quasi-smooth.
\end{proposition}
\begin{proof}
We first prove that $\ALG$ is 2-homogeneous, taking a square-zero extension $A \to B$ and a morphism $C \to B$. Homogeneity of $\HHom_{\ALG}$ is almost immediate, because $S\in \ALG(A\by_BC)$ is flat over $A\by_BC$ as a graded algebra (resp. levelwise flat over $A\by_BC$), giving the pullback square
\[
\begin{CD}
\HHom_{\ALG(A\by_BC)}(S,S') @>>> \HHom_{\ALG(A)}(S\ten_{A\by_BC}A ,S'\ten_{A\by_BC}A) \\
@VVV @VVV\\
\HHom_{\ALG(C)}(S\ten_{A\by_BC}C,S'\ten_{A\by_BC}C) @>>> \HHom_{\ALG(B)}(S\ten_{A\by_BC}B,S'\ten_{A\by_BC}B) .
\end{CD}
\]
Given objects $S_A \in \ALG(A)$ and $S_C \in \ALG(C)$, and an isomorphism $\alpha:S_A\ten_AB \to S_C\ten_CB$, we can define $S \in \ALG(A\by_BC)$ by $S:= S_A\by_{\alpha, S_C\ten_CB}S_C$, noting that this is automatically cofibrant, with $S\ten_{A\by_BC}A \cong S_A$ and $S\ten_{A\by_BC}C \cong S_C$. This completes the proof that $\ALG$ is 2-homogeneous.

Now, since any $S \in \ALG(A)$ is cofibrant over $A$, the map $f:\HHom_{\ALG(A)}(S, S') \to \HHom_{\ALG(B)}(S\ten_AB, S'\ten_AB)$ is a fibration for all square-zero extensions $A \to B$, because $S' \to S'\ten_AB$ is also a square-zero extension. If $A \to B$ is also a weak equivalence, then so is $S' \to S'\ten_AB $ (as $S'$ is levelwise flat over $A$, so the kernel is $S'\ten_AI$), and $f$ is therefore a trivial fibration. We have therefore shown that $\HHom_{\ALG}(S,S')$ is formally quasi-smooth

It only remains to show that $\ALG(A) \to \ALG(B)$ satisfies Condition (F2) from Definition \ref{2fibrn} for all square-zero extensions $A \to B$, and satisfies condition (TF2) whenever the extension is acyclic. 
\cite{stacks2} Proposition \ref{stacks-deform1} shows that the obstruction (in the simplicial case) to lifting $S \in \ALG(B)$ to $\ALG(A)$ lies  in 
\[
\EExt^2_B(\Omega^{S/B},S\ten_BI)=  \DD^2_{\id_S}(\HHom_{\ALG}(S,S), I), 
\]
where $I= \ker(A \to B)$. The idea is that as an almost simplicial algebra (i.e.  discarding $\pd_0$),  $S$ has a unique levelwise flat lift to $\tilde{S}$ over $A$. We can then choose a lifting $\delta$ of $\pd_0$ with $\pd_i\delta= \delta\pd_{i-1}$ for all $i>1$, and the obstruction comes from comparing $\delta^2$ and $\delta\pd_1$. 
The same proof adapts to the dg case, where we lift $S$ as a flat graded algebra, lift $d$ to $\delta$, and get $\delta^2$ as the obstruction. In particular, this ensures that (TF2) is satisfied, since if $\H_*(I)=0$, then $\EExt^2_B(\Omega^{S/B},S\ten_BI)=0$.

We now establish (F2): take $S \in \ALG(A)$, write $\bar{S}:= S\ten_AB$, and take   a homotopy equivalence $\theta:\bar{S}  \to T$. We now proceed with the proof in the dg case --- the simplicial case is entirely similar, replacing graded algebras with almost simplicial algebras.
As a graded algebra,  $T$ has a unique lift  $\tilde{T}$ over $A$, and we can then lift $\theta$  to $\tilde{\theta}:S \to \tilde{T}$ (as a graded morphism), since $S$ is smooth. We also lift the differential $d$ on $T$ to $\delta$ on $\tilde{T}$. 
The obstruction to these giving a morphism of dg algebras is then 
\[
(u,v):= (\delta^2, \delta \circ  \tilde{\theta} - \tilde{\theta} \circ d) \in \HOM^2_B(\Omega^{T/B}, T\ten_BI) \by \HOM^1_B(\Omega^{\bar{S}/B}, T\ten_BI),
\] 
which satisfies $du=0, dv= u \circ \theta$. Here, $\HOM^i_B$ denotes $B$-linear graded morphisms of chain degree $-i$.
Alternative choices for $(\delta, \tilde{\theta})$ would be of the form $ (\delta+a, \tilde{\theta}+b)$, for 
\[
(a,b) \in \HOM^1_B(\Omega^{T/B}, T\ten_BI) \by \HOM^0_B(\Omega^{\bar{S}/B}, T\ten_BI);
\]
these send $(u,v)$ to $(u+da, v+db+ a\circ \theta)$. Thus the obstruction to lifting $\theta$ and $T$ compatibly lies in $\H^2$ of the cone complex 
\[
 \HOM^{\bt}_B(\Omega^{T/B}, T\ten_BI)\xra{\theta^*} \HOM^{\bt}_B(\Omega^{\bar{S}/B}, T\ten_BI).
\]
Since $\theta$ is a weak equivalence, $\theta^*$ is a quasi-isomorphism, so cohomology of this complex is $0$, and suitable liftings thus exist.
\end{proof}

\begin{definition}\label{opencat}
Given a functor $\C: \Alg_{\H_0R} \to s\Cat$ and a functorial simplicial subcategory $\cM \subset \C$, say that $\cM$ is an \emph{open} simplicial subcategory of $\C$ if
\begin{enumerate}
\item $\cM$ is a full simplicial subcategory, i.e. for all $X, Y \in \cM(A)$, the map $\HHom_{\cM(A)}(X,Y) \to \HHom_{\C(A)}(X,Y)$ is a weak equivalence;

\item  the map $\cM \to \C$ is homotopy formally \'etale, in the sense that for any square-zero extension $A \to B$, the map
\[
\pi_0\cM(A) \to \pi_0\C(A)\by^{(2)}_{\pi_0\C(B)}\pi_0\cM(B)
\]
is essentially surjective on objects.
\end{enumerate}
\end{definition}

\begin{definition}
Given a  simplicial category $\C$ for which $\pi_0\C$ is a groupoid, say that $\C$ is  $n$-truncated if 
\[
\pi_i\HHom_{\C(A)}(X,Y)=0
\]
for all $i\ge n$. Thus $\C$ is $n$-truncated if and only if $\pi_i(\bar{W}\C)=0 $ for all $i>n$.
\end{definition}

\begin{corollary}\label{representaffine}
Take an  $n$-truncated presheaf $\cM: \Alg_{\H_0R} \to s\Cat$, open in $\pi^0\cW(\ALG) $,
such that 
\begin{enumerate}
\item[(0)] if $\{f_i \co A \to B_i\}_{i \in I}$ is an \'etale cover in $\Alg_{H_0R}$, then an object   $X \in \pi_0\cW(\ALG)(A)$ lies in the essential image of $\pi_0\cM(A)$ whenever $f_i(X)$ lies in the essential image of $\pi_0\cM(B_i)$ for all $i$.

\item
for all finitely generated $A \in \Alg_{\H_0R}$  and all $X \in \cM(A)$, the functors 
\[
 \EExt^i_{X}(\bL^{X/A},X\ten_A^{\oL}- ): \Mod_A \to \Ab
\]
 preserve filtered colimits for all $i>1$.

\item for all finitely generated integral domains $A \in \Alg_{\H_0R}$  and all $X \in \cM(A)$, the groups $\EExt^i_{X}(\bL^{X/A}, X)$ are all  finitely generated $A$-modules.

\item for all complete discrete local Noetherian  $\H_0R$-algebras $A$, with maximal ideal $\m$, the map
$$
\cM(A) \to {\Lim}^h \cM(A/\m^n)
$$
is a weak equivalence (see  Remark \ref{formalexistrk} for a reformulation).

\item
The functor $c(\pi_0\cM):  \Alg_{\pi_0R} \to \Set$  of components of the groupoid $\pi_0\cM$ preserves filtered colimits.

\item
For all $A \in \Alg_{\pi_0R}$ and all $X \in \cM(A)$, the functors $\pi_i\HHom_{\cM}(X,X): \Alg_A \to \Set$  preserve filtered colimits for all $i\ge 0$.
\end{enumerate}

Let $\tilde{\cM}:d\cN^{\flat}_R \to s\Cat$  be the full simplicial subcategory of $\cW(\ALG)(A)$ consisting of objects $C$ for which $C\ten_A\H_0A$ is weakly equivalent (in  $\cW(\ALG)(\H_0A)$) to an object of  $\cM(\H_0A)$. 
Then the nerve $\bar{W}\tilde{\cM}$ is a geometric derived $n$-stack.
\end{corollary}
\begin{proof}
The proof of \cite[Lemma \ref{stacks-dgmodhyp}]{stacks2} adapts from modules to algebras,  showing that $A \mapsto dg_+\Alg_A$ and $A \mapsto s\Alg_A$ are left Quillen hypersheaves, and hence that $\bar{W}\cW\ALG$ is a hypersheaf by \cite[Proposition \ref{stacks-sheafprop}]{stacks2}. Then Proposition \ref{descentlemma} ensures that $\bar{W}\cW\cM$ is an \'etale hypersheaf.

Since $\pi_0\cW(\ALG)$ is a groupoid-valued functor, and $\cM \into \pi^0\cW(\ALG)$ is open,  the map
\[
 \cM \to \pi^0\cW(\ALG)
\] 
is homotopy formally \'etale. This also implies that $\bar{W}\tilde{\cM}$ is a model for
\[
A\mapsto \bar{W}\cM(\H_0A)\by^h_{\bar{W}\cW(\ALG)(\H_0A)}\bar{W}\cW(\ALG)(A).
\]

Now observe that this satisfies the conditions of Theorem \ref{lurierep3}, using Lemmas \ref{2htpicgood} and \ref{2hgsgood}, Propositions \ref{wetale} and \ref{alghgs}, and Corollary \ref{w2good}. 
The proof of Proposition \ref{alghgs} shows that 
\[
\DD^i_X(\bar{W}\tilde{\cM}/R, M)= \EExt^{i+1}_{X}(\bL^{X/A},X\ten_A^{\oL}M ).
\]
\end{proof}

We now give a simple motivating example. This is already well-known; indeed, \cite{hag2} Proposition 2.2.6.8 generalises it to algebras for any operad. 

\begin{example}[Moduli of finite schemes]\label{fineg}
Define $\cM: \Alg_{\H_0R} \to s\Cat$ by setting
 $\cM(A) \subset \pi^0\cW(\ALG)(A)$ to consist of dg (resp. simplicial) $A$-algebras $S$ for which $\H_0S$ is a finite flat $A$-module, and $\H_iS=0$ for all $i>0$. Then  we get
\[
\tilde{\cM}\simeq \cM\by^{(2)}_{\pi^0\ALG}\ALG.
\]
Now, $\H_0$ gives a functor from $\cM(A)$ to the groupoid $\cW(\FALG(A))$ of finite flat $A$-algebras, and this is a weak equivalence in $s\Cat$, so 
\[
\cM \simeq \cW(\FALG),
\] 
which is a hypersheaf by \cite{drep} Remark \ref{drep-stackhyper}. Thus the corollary above implies that the nerve $\bar{W}\tilde{\cM}$
is a geometric derived $1$-stack, which 
we  can now regard as the derived moduli stack of finite flat algebras.
\end{example}

We now check  that this example is consistent with the flatness criteria used in \cite{hag2} and \cite{lurie}.

\begin{definition}
Given $A \in dg_+\Alg_R$ (resp. $A \in s\Alg_R$) say that $M\in dg_+\Mod_A$  (resp. $M \in s\Mod_A$) is homotopy-flat over $A$ if $ \H_0M$ is flat over $\H_0A$, and the maps
\[
\H_n(A)\ten_{\H_0A}\H_0M \to \H_nM
\]
are isomorphisms for all $n$.
\end{definition}

Note that \cite{hag2} Lemma 2.2.2.2 
implies that this is the same as flatness in the sense of \cite{hag2} Definition 1.2.4.1, which corresponds to the functor $\ten_A^{\oL}M: \Ho(dg_+\Mod_A) \to \Ho(dg_+\Mod_A)$ (resp. $\ten_A^{\oL}M: \Ho(s\Mod_A) \to \Ho(s\Mod_A)$) preserving homotopy pullbacks. 

\begin{lemma}
A morphism $f:A \to B$ in  $A \in dg_+\Alg_R$ (resp. $s\Alg_R$) is homotopy-flat if and only if
$B\ten_A^{\oL}\H_0A$ is weakly equivalent to a discrete flat $\H_0A$-algebra.
\end{lemma}
\begin{proof}
If $f$ is flat, then the proof of \cite{hag2} Lemma 2.2.2.2(2) shows that for any $A$-module $N$, we have $\H_*(\oL f^*N)\cong \H_0B\ten_{\H_0A}\H_*N$. Setting $N=\H_0A$, this gives $ \oL f^*\H_0A \simeq \H_0B$. This is flat over $\H_0A$ by hypothesis, completing the proof of one implication. 

For the converse, assume that $ \oL f^*\H_0A \simeq \H_0B$, flat over $\H_0A$. By considering  free simplicial resolutions, this implies that $\oL f^*M \simeq \H_0B\ten_{\H_0A}M$ for all $\H_0A$-modules $M$. We have a spectral sequence
\[
\H_i(\oL f^* \H_jA) \abuts \H_{i+j}(B),
\]
giving $\H_jB \cong \H_0B\ten_{\H_0A}\H_jA$, which completes the proof.
\end{proof}

In particular, this means that $B$ is homotopy-flat over $A$ if and only if for all $\H_0A$-algebras $C$ and all $i>0$, we have $\H_i(B\ten_A^{\oL}C)=0$. 

\begin{proposition}\label{relalghgs}
For $T \in dg_+\Alg_R$ (resp. $T \in s\Alg_R$) flat as a graded algebra over $R$, the functor $\ALG_T:=\ALG(T\ten -):  d\cN^{\flat} \to s\Cat$  is 2-homogeneous and formally 2-quasi-smooth.
\end{proposition}
\begin{proof}
The proof of Proposition \ref{alghgs} carries over to this generality, noting that for any (acyclic) square-zero extension $A \to B$, we get another (acyclic) square-zero extension $T\ten A \to T\ten B$. 
\end{proof}

\begin{corollary}\label{representrelaffine}
Take an  $n$-truncated presheaf $\cM: \Alg_{\H_0R} \to s\Cat$, open in $\pi^0\cW(\ALG_T) $,
such that 
\begin{enumerate}
\item[(0)] if $\{f_i \co A \to B_i\}_{i \in I}$ is an \'etale cover in $\Alg_{H_0R}$, then an object   $X \in \pi_0\cW(\ALG_T)(A)$ lies in the essential image of $\pi_0\cM(A)$ whenever $f_i(X)$ lies in the essential image of $\pi_0\cM(B_i)$ for all $i$.

\item
for all finitely generated $A \in \Alg_{\H_0R}$  and all $X \in \cM(A)$, the functors 
\[
 \EExt^i_{X}(\bL^{X/A\ten T},X\ten_A^{\oL}- ) :\Mod_A \to \Ab
\]
 preserve filtered colimits for all $i>0$.

\item for all finitely generated integral domains $A \in \Alg_{\H_0R}$  and all $X \in \cM(A)$, the groups $\EExt^i_{X}(\bL^{X/A\ten T}, X)$ are all  finitely generated $A$-modules.

\item for all complete discrete local Noetherian  $\H_0R$-algebras $A$, with maximal ideal $\m$, the map
$$
\cM(A) \to {\Lim}^h \cM(A/\m^n)
$$
is a weak equivalence (see  Remark \ref{formalexistrk} for a reformulation).

\item
The functor $c(\pi_0\cM):  \Alg_{\pi_0R} \to \Set$  of components of the groupoid $\pi_0\cM$ preserves filtered colimits.

\item
For all $A \in \Alg_{\pi_0R}$ and all $X \in \cM(A)$, the functors $\pi_i\HHom_{\cM}(X,X): \Alg_A \to \Set$  preserve filtered colimits for all $i\ge 0$.

\end{enumerate}

Let $\tilde{\cM}:d\cN^{\flat}_R \to s\Cat$  be the full simplicial subcategory of $\cW(\ALG_T)(A)$ consisting of objects $C$ for which $C\ten_A\H_0A$ is weakly equivalent (in  $\cW(\ALG_T)(\H_0A)$) to an object of  $\cM(\H_0A)$. 
Then the nerve $\bar{W}\cM$ is a geometric derived $n$-stack.
\end{corollary}
\begin{proof}
The proof of Corollary \ref{representaffine} carries over.
\end{proof}

\begin{example}[Moduli of finite schemes over an affine base]\label{felfineg}
Fix $T \in dg_+\Alg_R$ (resp. $T \in s\Alg_R$) flat as a graded algebra over $R$, define $\cM: \Alg_{\H_0R} \to s\Cat$ by
setting $\cM(A) \subset \pi^0\cW\ALG(T\ten_R A)$ to consist of dg (resp. simplicial) $T\ten_RA$-algebras $S$ for which $\H_0S$ is a finite flat $A$-module, and $\H_iS=0$ for all $i>0$. Then  we get
\[
\tilde{\cM} \simeq  \cM\by^{(2)}_{\pi^0\ALG_T}\ALG_T.
\]
Now, $\H_0$ gives a functor from $\cM(A)$ to the groupoid $\cW( \FALG_{\H_0T}(A))$ of $(\H_0T\ten_{\H_0R}A)$-algebras, finite and flat over $A$, and this is a weak equivalence in $s\Cat$, so 
\[
\cM \simeq \cW\FALG_{\H_0T},
\] 
which is a hypersheaf by \cite{drep} Remark \ref{drep-stackhyper}. Thus the corollary above implies that the nerve
$
\bar{W}\tilde{\cM}
$
is a geometric derived $1$-stack, which 
we  can now regard as the derived moduli stack of finite flat algebras equipped with a morphism from $T$.
\end{example}

\subsection{Derived moduli of schemes}

Recall that we write $d\Alg_R$ for either $dg_+\Alg_R$ or $s\Alg_R$.

\subsubsection{Cosimplicial derived rings}

\begin{definition}
 Define $\Delta$ to be the  ordinal number category. For any category $\C$, this means that $\C^{\Delta}$ (resp. $\C^{\Delta^{\op}}$) is the category of cosimplicial (resp. simplicial) diagrams in $\C$.
\end{definition}

\begin{definition}
Define a simplicial structure on  $(d\Alg_A)^{\Delta}$ as follows.  For $C \in(d\Alg_A)^{\Delta}$,  first  define $h(\Delta^n,C) \in (d\Alg_A)^{\Delta}$
by $h(\Delta^n,C)^i:= (C^i)^{\Delta^n_i}$, and then define the  
$\Hom$-spaces  by
\[
\HHom_{ (d\Alg_A)^{\Delta}}(B,C)_n =  \Hom_{(d\Alg_A)^{\Delta} }(B, h(\Delta^n,C)). 
\] 
\end{definition}

\begin{definition}
Given $A \in (d\Alg_R)^{\Delta} $, define $L^nA$ to represent the functor $M_n\Spec A$ 
given by
$$
(M_n\Spec A)(B) = \{ x \in \prod_{i=0}^n \Hom_{d\Alg_R}(A^{n-1},B)  \,:\, (\pd^i)^*x_j= (\pd^{j-1})^*x_i 
\text{ if } 0\le i<j\le n \}.
$$
\end{definition}

\begin{definition}
 Say that a morphism $A \to B$ in $(d\Alg_R)^{\Delta} $ is a Reedy cofibration  if the latching maps
\[
A^n\ten_{L^nA}L^nB \to B^n 
\]
in $d\Alg_A $  are cofibrations  (in the sense of Definitions \ref{dgcof} and \ref{scof}) for all $n \ge 0$.

This is equivalent to satisfying the  conditions that for  all acyclic surjections in $C' \to C$ in $d\Alg_R$, the map 
\[
\Hom_{d\Alg_R}(B,C') \to \Hom_{d\Alg_R}(B,C)\by_{\Hom_{d\Alg_R}(A,C)}\Hom_{d\Alg_R}(A,C')
\] 
is a trivial fibration in $\bS$, where we set $\Hom_{d\Alg_R}(A,C)_n:= \Hom_{d\Alg_R}(A^n,C)$.
\end{definition}

\begin{definition}
For  $A \in d\Alg_R$, 
define $c\ALG_T(A)$ to be the full simplicial subcategory of the category  $ (d\Alg_{(T\ten_RA) })^{\Delta} $  on cofibrant  objects (over $T\ten_RA$). 
\end{definition}

\begin{proposition}\label{calghgs}
For $T$ levelwise flat as a graded algebra over $R$, the functor $c\ALG_T:  d\cN^{\flat} \to s\Cat$ defined above is 2-homogeneous and formally 2-quasi-smooth.
\end{proposition}
\begin{proof}
This is a fairly direct consequence of the proof of Proposition \ref{relalghgs}. Again the only part which is not immediate is showing that $c\ALG_T(A) \to c\ALG_T(B)$ satisfies Condition (F2) from Definition \ref{2fibrn} for all square-zero extensions $A \to B$, and satisfies condition (TF2) whenever the extension is acyclic.

Given a square-zero  extension $A \to B$ and an object $S$  in $c\ALG_T(B)$,     Proposition \ref{relalghgs} gives criteria for lifting $S^0$ to $\tilde{S}^0 \in \ALG_{T^0}(A)$. Similar criteria apply for lifting the $L^0S\ten_{L^0T}T^1$-algebra $S^1$ to an $L^0\tilde{S}\ten_{L^0T}T^1$-algebra $\tilde{S}^1$, with an augmentation $\tilde{S}^1 \to M^1(\tilde{S}) $ extending the augmentation $S^1 \to  M^1S$. 
 We then proceed by induction on the levels of the cosimplicial diagram. A similar argument applies to lift morphisms.
\end{proof}

\begin{definition}\label{dsmooth}
Say that a map $f:A \to B$ in $d\Alg_R $  is  homotopy-flat (resp. homotopy-smooth, resp. homotopy-\'etale)  
 if $\H_0A \to \H_0B$ is flat (resp. smooth, resp. \'etale)  and the maps
\[
\H_nA\ten_{\H_0A} \H_0B \to \H_nB 
\]
are isomorphisms for all $n$.

Say that a homotopy-flat map $A \to B$ is a cover if  $\H_0A \to \H_0B$ is  faithfully flat.
\end{definition}

\begin{definition}\label{cofetale}
Say that a Reedy cofibration $A \to B$ in $(d\Alg_R)^{\Delta} $  is a homotopy-\'etale hypercover if  the latching maps
\[
A^n\ten_{L^nA}L^nB \to B^n 
\]
in $d\Alg_A $  are homotopy-\'etale covers for all $n \ge 0$.
\end{definition}

By \cite{stacks2} Proposition \ref{stacks-procofibrant}, there is a functor $Q$ from $c\Alg_T(A)$ to the category $\ind( c\Alg_T(A))$ of direct systems, equipped with a natural transformation $\id \to Q$. This has the properties that $QC= \{(QC)_{\alpha}\}_{\alpha}$ is a system of  Reedy cofibrant homotopy-\'etale hypercovers of $C $, and that any such hypercover of $QC$ has a retraction.

\begin{definition}
Define $c\ALG_{T,\hyp}:d\cN^{\flat} \to s\Cat$ by setting  $c\ALG_{T,\hyp}(A)$ to have the same objects as $c\Alg_T(A)$, and morphisms
\begin{align*}
\HHom_{c\ALG_{T,\hyp}(A)}(B,C) &=  \HHom_{\ind(c\Alg_T(A))}(QB, QC)\\
&= \Lim_{\beta} \LLim_{\alpha} \HHom_{\ind(c\Alg_T(A))}((QB)_{\beta}, (QC)_{\alpha})\\
&\simeq \LLim_{\alpha} \HHom_{\ind(c\Alg_T(A))}(B, (QC)_{\alpha}).
\end{align*}
\end{definition}

\begin{proposition}\label{calghyphgs}
The functor $c\ALG_{T,\hyp}:  d\cN^{\flat} \to s\Cat$ defined above is 2-homogeneous and formally 2-quasi-smooth. 
\end{proposition}
\begin{proof}
This follows immediately from Proposition \ref{calghgs}, making use of the fact that filtered colimits commute with finite limits.
\end{proof}

\subsubsection{Artin hypergroupoids}

\begin{definition}
Define the category $d\Aff_R$ of derived affine $R$-schemes to be opposite to  $d\Alg_R$ (i.e. $dg_+\Alg_R$ or $s\Alg_R$). Write $sd\Aff_R$ for the category of simplicial objects in $d\Aff_R$; this is opposite to $(d\Alg_R)^{\Delta}$. Say that a morphism in $d\Aff_R$ is a fibration if the corresponding morphism in $d\Alg_R$  is a cofibration.
\end{definition}

The following definition is adapted from \cite{stacks2}. For the definitions of $\pd\Delta^m$ and $\L^m_k$, see Definition \ref{horndef}.
\begin{definition}\label{dnpreldef}
Given   $Y_{\bt} \in sd\Aff_R$, define a     derived Artin (resp. Deligne--Mumford) $n$-hypergroupoid  over $Y_{\bt}$ to be a morphism $X_{\bt}\to Y_{\bt}$ in $sd\Aff_R$, satisfying the following conditions:
\begin{enumerate}

\item for each  $m\ge 0$, the  matching map
$$
X_m \to \Hom_{\bS}(\pd\Delta^m, X)\by_{\Hom_{\bS}(\pd\Delta^m,Y)}Y_m 
$$
is a fibration  (this is equivalent to saying that $X \to Y$ is opposite to a Reedy cofibration in   $(d\Alg_R)^{\Delta}$);

\item for each  $k,m$, the  partial matching map
$$
X_m \to \Hom_{\bS}(\L^m_k, X)\by_{\Hom_{\bS}(\L^m_k,Y)}Y_m 
$$
is a homotopy-smooth (resp. homotopy-\'etale) cover (in the sense of Definition \ref{dsmooth});

\item for all $m>n$ and all $k$, the  partial matching maps
$$
X_m \to \Hom_{\bS}(\L^m_k, X)\by_{\Hom_{\bS}(\L^m_k,Y)}Y_m 
$$
are weak equivalences.
\end{enumerate}
\end{definition}

\begin{proposition}\label{bigthmd}
The simplicial category  of strongly quasi-compact $n$-geometric derived Artin (resp. Deligne--Mumford) stacks over the hypersheafification $\Spec (T\ten_R^{\oL} A)^{\sharp}$ is weakly equivalent to the opposite of the full subcategory $\HYP^n_T(A)$ of $c\ALG_{T,\hyp}(A)$ whose objects are dual to   derived Artin (resp. Deligne--Mumford) $n$-hypergroupoids. 
\end{proposition}
\begin{proof}
This is \cite{stacks2} Theorem \ref{stacks-bigthm}. The functor is given by sending $C \in \HYP^n(A)$ to the hypersheafification of the presheaf $\Spec C: d\Aff_R^{\op} \to \bS$.
\end{proof}

\begin{definition}
Given a flat $\infty$-geometric  Artin stack $\fY$ over $\H_0R$, and $A \in \Alg_{\H_0R}$,  let $F\Stack^n_{\fY}(A)$ be the simplicial category of  $n$-geometric  Artin stacks over $\fY\ten_{\H_0R}A$, flat over $A$,  and set $F\Stack_{\fY}(A)= \LLim_n F\Stack^n_{\fY}(A)$. Likewise,
given an $\infty$-geometric derived Artin stack $\fY$ over $R$ and $A \in d\Alg_R$, let $d\Stack^n_{\fY}(A)$ be the simplicial category of $n$-geometric derived Artin stacks over $\fY\ten_R^{\oL}A$, and set $d\Stack_{\fY}(A)= \LLim_n d\Stack^n_{\fY}(A)$.
\end{definition}

Note that if we set $\fY:= \Spec (T)^{\sharp}$, then Proposition \ref{bigthmd} shows that sheafification gives an $\infty$-equivalence $\HYP^n_T(A) \to d\Stack^n_{\fY}(A)$. 
Also note that the maps  $d\Stack^n_{\fY}(A) \to d\Stack_{\fY}(A)$ are homotopy formally \'etale, as are the maps   $F\Stack^n_{\fY}(A) \to d\Stack^n_{\fY}(A)$,   $F\Stack_{\fY}(A) \to d\Stack_{\fY}(A)$
and the embedding $F\Sch_{\fY}(A) \to F\Stack_{\fY}(A)$ of algebraic spaces into  $\infty$-geometric stacks.

\begin{lemma}\label{dstack}
Given a strongly quasi-compact $\infty$-geometric derived Artin stack $\fY$ over $R$, the functor $ d\Stack_{\fY}:d\cN^{\flat}_R \to s\Cat $ is homotopy-preserving and homotopy-homogeneous.
\end{lemma}
\begin{proof}
Choose a representative $T \in c\ALG(R)$ with $(\Spec T)^{\sharp} \simeq \fY$ (as in Proposition \ref{bigthmd}). If we write $\HYP_T(A):= \LLim_n \HYP^n_T(A)$, then  Proposition \ref{bigthmd} allows us to replace  $ d\Stack_{\fY}(A) $ with $\HYP_T(A)$. The map $\HYP_T \to c\ALG_{T,\hyp}(A)$ is homotopy formally \'etale, so it suffices to show that $c\ALG_{T,\hyp}(A)$ is homotopy-preserving and homotopy-homogeneous. This follows from combining Lemmas \ref{2htpicgood} and \ref{2hgsgood} with Proposition \ref{calghyphgs}.
\end{proof}

\begin{theorem}\label{representdaffine}
Take an  $\infty$-geometric  derived Artin stack $\fY$,  homotopy-flat over $R$. 

Assume that we have  an  $n$-truncated presheaf $\cM: \Alg_{\H_0R} \to s\Cat$, open in $ \cW( F\Stack_{\fY} )$ (in the sense of Definition \ref{opencat}),
noting that homotopy-flatness ensures that $Y_A:=\fY\ten^{\oL}_RA $ is a flat  $\infty$-geometric  (underived) Artin stack over $A$  for all $A \in \Alg_{\H_0R}$.

Also assume that this satisfies the following conditions
\begin{enumerate}
\item[(0)] if $\{f_i \co \Spec B_i \to \Spec A\}_{i \in I}$ is an \'etale cover in affine $\H_0R$-schemes, then an object  $X \in \pi_0\cW(F\Stack_{\fY})(A)$ lies in the essential image of $\pi_0\cM(A)$ whenever $f_i^*(X)$ lies in the essential image of $\pi_0\cM(B_i)$ for all $i$.

\item
for all finitely generated $A \in \Alg_{\H_0R}$  and all $X \in \cM(A)$, the functors 
\[
 \EExt^i_{X}(\bL^{X/Y_A},\O_{X}\ten_A^{\oL}- ): \Mod_A \to \Ab
\]
preserve filtered colimits for all $i>1$.

\item for all finitely generated integral domains $A \in \Alg_{\H_0R}$  and all $X \in \cM(A)$, the groups $\EExt^i_{X}(\bL^{X/Y_A}, \O_{X} )$ are all  finitely generated $A$-modules.

\item for all complete discrete local Noetherian  $\H_0R$-algebras $A$, with maximal ideal $\m$, the map
$$
\cM(A) \to {\Lim}^h \cM(A/\m^n)
$$
is a weak equivalence (see  Remark  \ref{formalexistrk} for a reformulation).

\item
The functor $c(\pi_0\cM):  \Alg_{\H_0R} \to \Set$  of components of the groupoid $\pi_0\cM$ preserves filtered colimits.

\item
For all $A \in \Alg_{\H_0R}$ and all $X \in \cM(A)$, the functors $\pi_i\HHom_{\cM}(X,X): \Alg_A \to \Set$  preserve filtered colimits for all $i\ge 0$.

\end{enumerate}

Let $\tilde{\cM}:d\cN^{\flat}_R \to s\Cat$  be the full simplicial subcategory of $\cW(c\ALG_{T,\hyp}(A))$ consisting of objects $\fX$ for which $\fX\ten_A\H_0A$ is weakly equivalent  in $d\Stack_{\fY}(\H_0A)$  to an object of  $\cM(\H_0A)$. 
Then the nerve $\bar{W}\tilde{\cM}$ is a geometric derived $n$-stack.
\end{theorem}
\begin{proof}
As in Lemma \ref{dstack}, we may replace $c\ALG_{T,\hyp}(A)$ with $ d\Stack_{\fY}(A) $, so
\[
\tilde{\cM}(A) \simeq \cM(\H_0A)\by^h_{\cW( d\Stack_{\fY}(\H_0A))}\cW( d\Stack_{\fY}(A)),
\] 
because $\cM $ is open in $F\Stack_{\fY} $, which is open in $\pi^0d\Stack_{\fY}$. 

The proof of Corollary \ref{representaffine} now  carries over to this context, using Lemma \ref{dstack} in place of Proposition \ref{alghgs}.  The only non-trivial part is the calculation 
\[
\DD^i_X(\bar{W}\tilde{\cM}/R, M) =\EExt^{i+1}_X(\bL^{X/Y_A}, \O_X\ten_A^{\oL}M),
\]
 but this follows from the calculations of \cite{stacks2} \S \ref{stacks-cotsn}, for instance.
\end{proof}

\begin{remark}\label{rep2rk}
Since $\bar{W}\tilde{\cM}$ is a geometric derived stack, it extends canonically to a functor
\[
\bar{W}\tilde{\cM}: d\Alg_R \to \bS
\]
where $d\Alg_R= s\Alg_R$ or $dg_+\Alg_R$. The obvious candidate for this functor is given by adapting the description above;  $\tilde{\cM}$  should be the  full simplicial subcategory of $\cW(c\ALG_{T,\hyp}) $ (or equivalently of $\cW( d\Stack_{\fY})$) consisting of objects $\fX$ for which $\fX\ten_A^{\oL}\H_0A$ is weakly equivalent  in $d\Stack_{\fY}(\H_0A)$  to an object of  $\cM(\H_0A)$.

To check that this is indeed the case, it would suffice
 to verify the conditions of \cite{drep} Corollary \ref{drep-lurierep2} for this functor. These reduce to showing that the functor
\[
\bar{W}\cW( d\Stack_{\fY}) : d\Alg_R \to \bS
\]
is  a homotopy-preserving functor (which is certainly true), and
nilcomplete in the sense that for all $ A \in d\Alg_R$, the map
$$
\bar{W}\cW( d\Stack_{\fY})(A) \to {\Lim}^h \bar{W}\cW( d\Stack_{\fY})(P_kA)
$$
is an equivalence, for $\{P_kA\}$ the Postnikov tower of $A$. Nilcompleteness certainly seems plausible, as similar statements in \cite{lurie}   \S 8 are asserted to be automatic. However, broadly similar comparisons in \cite{hag2} involve weighty strictification theorems (such as \cite{hag2} Theorem B.0.7 or \cite{TVsegal} Theorem 4.2.1). 
\end{remark}

We now give a generalisation which will be used for moduli of pointed schemes.

\begin{definition}
Given a morphism $ f: \fZ \to \fY$ of  $\infty$-geometric derived Artin stacks over $R$,  define $d\Stack_{\fY}^{\fZ}:d\cN^{\flat}_R \to s\Cat  $ as follows.  First set $\fZ_A:= \fZ\ten_R^{\oL}A$, $\fY_A:= \fY\ten_R^{\oL}A$;  objects of  $d\Stack_{\fY}^{\fZ}(A)$ are factorisations 
\[
\fZ_A \xra{i} X \xra{p}   \fY_A
\] 
of $f$ in the category of $\infty$-geometric derived Artin stacks over $A$. Morphism spaces are given by homotopy fibre products 
\[
\HHom_{d\Stack_{\fY}^{\fZ}(A)}(X, X'):= \HHom_{\fY_A}(X, X') \by^h_{i^*,\HHom_{\fY_A}(\fZ_A , X')} \{i'\}.
\]
\end{definition}

\begin{theorem}\label{augrepresentdaffine}
Take a morphism $\fZ \to \fY$ of   homotopy-flat $\infty$-geometric  derived Artin stacks   over $R$. 
Assume that we have  an  $n$-truncated presheaf $\cM: \Alg_{\H_0R} \to s\Cat$, open in the functor $A \mapsto  \cW((\fZ\ten^{\oL}_RA)\da  F\Stack_{\fY} )$ (in the sense of Definition \ref{opencat}),
noting that  $Y_A:=\fY\ten^{\oL}_RA $ and $Z_A:=\fZ\ten^{\oL}_RA $ are flat  $\infty$-geometric  (underived) Artin stacks over $A$  for all $A \in \Alg_{\H_0R}$.

Also assume that these  satisfy the following conditions
\begin{enumerate}
\item[(0)] if $\{f_i \co \Spec B_i \to \Spec A\}_{i \in I}$ is an \'etale cover in affine $\H_0R$-schemes, then an object  $X \in \pi_0\cW((\fZ\ten^{\oL}_RA)\da  F\Stack_{\fY} )(A)$ lies in the essential image of $\pi_0\cM(A)$ whenever $f_i^*(X)$ lies in the essential image of $\pi_0\cM(B_i)$ for all $i$.

\item
for all finitely generated $A \in \Alg_{\H_0R}$  and all $ Z_A\xra{i} X$ in $\cM(A)$, the functors 
\[
 \EExt^j_{X}(\bL^{X/Y_A},(\O_{X}  \to \oR i_*\O_{Z_A})\ten_A^{\oL}- ): \Mod_A \to \Ab
\]
 preserve filtered colimits for all $j>1$.

\item for all finitely generated integral domains $A \in \Alg_{\H_0R}$  and all $X \in \cM(A)$, the groups $\EExt^j_{X}(\bL^{X/Y_A}, \O_{X} \to \oR i_*\O_{Z_A})$ are all  finitely generated $A$-modules.

\item for all complete discrete local Noetherian  $\H_0R$-algebras $A$, with maximal ideal $\m$, the map
$$
\cM(A) \to {\Lim}^h \cM(A/\m^n)
$$
is a weak equivalence (see Remark \ref{formalexistrk} for a reformulation).

\item
The functor $c(\pi_0\cM):  \Alg_{\pi_0R} \to \Set$  of components of the groupoid $\pi_0\cM$ preserves filtered colimits.

\item
For all $A \in \Alg_{\pi_0R}$ and all $X \in \cM(A)$, the functors $\pi_j\HHom_{\cM}(X,X): \Alg_A \to \Set$  preserve filtered colimits for all $j\ge 0$.
\end{enumerate}

Let $\tilde{\cM}:d\cN^{\flat}_R \to s\Cat$  be the full simplicial subcategory of $\cW(d\Stack_{\fY}^{\fZ} )$ consisting of morphisms $i:Z_A\to \fX$ for which $\fX\ten_A\H_0A$ is weakly equivalent  in $ d\Stack_{\fY}^{\fZ }(\H_0A)$  to an object of  $\cM(\H_0A)$. 
Then the nerve $\bar{W}\tilde{\cM}$ is a geometric derived $n$-stack.
\end{theorem}
\begin{proof}
Lemma \ref{dstack} ensures that $d\Stack_{\fY}$ is homotopy-preserving and homotopy-homogeneous, so it follows immediately that $d\Stack_{\fY}^{\fZ} $ inherits these properties. Lemma \ref{wetale} ensures that the same is true of $\cW(d\Stack_{\fY}^{\fZ} )$.

In order to verify the conditions of Theorem \ref{lurierep3}, we now just need to establish the isomorphisms
\[
\DD^j_{(Z_A \xra{i} X)}(\bar{W}\tilde{\cM}/R, M) =\EExt^{j+1}_X(\bL^{X/Y_A}, \O_X\ten_A^{\oL}M \to \oR i_* \O_{Z_A}\ten_R^{\oL}M).
\]
As in the proof of Theorem \ref{representdaffine}, we know that the loop space of $T_{(Z_A \xra{i} X)}(\bar{W}\tilde{\cM}/R)(M)$ is the homotopy fibre of
\[
\oR\HHom_X(\bL^{X/Y_A}, \O_X\ten_A^{\oL}M) \to \oR\HHom_{Z_A}(\oL i^*\bL^{X/Y_A}, \O_{Z_A}\ten_R^{\oL}M )
\]
over $0$. We can rewrite the second space as $\oR\HHom_{X}(\bL^{X/Y_A}, \oR i_*\O_{Z_A}\ten_R^{\oL}M )$, and the result follows by noting that homotopy fibres correspond to mapping cones of complexes. 
\end{proof}

\subsection{Examples}
In the examples which follow, we will always  choose a functorial full subgroupoid $\cM$ of the core $\cW(F\Stack_{Y}^n)$ on objects which are fibred in geometric $0$-stacks (i.e. algebraic spaces) over $Y\ten_{\H_0R}A$ --- this has the effect that $\cM $ is a bona fide groupoid, not just an $\infty$-groupoid.  We get a derived $1$-stack by taking the nerve $\bar{W}\tilde{\cM}$ of $\tilde{\cM}$ as constructed in Theorem \ref{representdaffine}.

Our first application generalises the examples of Hilbert schemes and moduli of curves in \cite{lurie} \S 8, but beware that the derived schemes featuring  in \cite{lurie} form a larger category than  derived geometric $0$-stacks.
\begin{example}[Moduli of schemes]\label{modscheg}
Let  $\fY$ be a homotopy-flat derived geometric  $0$-stack   over $R$, so $Y:=\fY\ten^{\oL}_R\H_0R$ is a flat algebraic space over $\H_0R$.  Choose a  stack $\cM$ equipped with a functorial fully faithful embedding   of $\cM$ into the core of  $F\Sch_{Y}$, with the property that $\cM$ is closed under infinitesimal deformations. This ensures that $\bar{W}\cM \to   \bar{W}F\Sch_{Y}$ is open.

An example satisfying these conditions is when $\fY=R$ and $\cM$ is some  moduli stack of schemes with fixed dimension. Another example is to let $\cM(A)$  be  the set of closed subschemes of $Y\ten_{\H_0R}A$. 

Then we can define $\tilde{\cM}$ as above, and we get a derived $1$-stack
$
\bar{W}\cM
$,
with $\pi^0\bar{W}\tilde{\cM} = B\cM$.

 If $\cM$ is known to be an algebraic stack, then the last two conditions of Theorem \ref{representdaffine} are automatically satisfied, so for geometricity of $\bar{W}\tilde{\cM}$ we need only check properties of the cohomology groups
\[
\DD^i_{f}( \bar{W}\tilde{\cM}/R, M)=\EExt^{i+1}_{X}(\bL^{X/Y\ten_{\H_0R} A}, \O_{X}\ten^{\oL}_AM ),
\]
for  $A \in \Alg_{\H_0R}$ and $X \in \cM(A)$.
\end{example}

Our next application gives a generalisation of \cite{hag2} Theorem 2.2.6.11. In particular, \cite{hag2} Corollary 2.2.6.14 addresses the case when $X$ is a flat projective scheme and $Y$ a smooth projective scheme. 

\begin{example}[Derived $\Hom$-stacks] \label{Mor}
Another example is to fix  $n$-geometric derived stacks $\fX,\fZ$, homotopy-flat over $R$. Set $X:= \pi^0\fX=\fX\ten^{\oL}_R\H_0R$, $Z:= \pi^0\fZ$; these are flat  $n$-geometric stacks over $\H_0R$.

Now, let $\fY:= \fX\by_{\Spec R} \fZ$ with $Y:= \pi^0\fY$,  and  set $\cM(A) \subset \Sch_Y(A)$ to consist of maps 
\[
U \to X \by_{\Spec \H_0R} Z\by_{\Spec \H_0R} \Spec A
\] 
for which the projection $U \to X\by_{\Spec \H_0R} \Spec A$ is an equivalence. Then $\bar{W}\cM(A)$ is isomorphic to the simplicial set $\HHom(X\by_{\Spec \H_0R} \Spec A,  Z)$. The derived moduli stack $\bar{W}\tilde{\cM}(A)$ is weakly equivalent to the nerve of
\[
d\Stack_{\fY}(A)\by^h_{d\Stack(A)}\{\fZ\}
\]
--- $\bar{W}\tilde{\cM}$ is the derived $\Hom$-stack $\cHom(\fX,\fZ)$.

By considering the arrow category of $d\Stack(A)$, we see that this in turn is equivalent to  $\HHom(\fX\by_{\Spec R} \Spec A, \fZ)$. Thus a far simpler characterisation of this moduli functor is given by choosing $B,C \in \HYP^n_R$ with $\fX= (\Spec B)^{\sharp},\fX= (\Spec C)^{\sharp}$. Then
\[
\tilde{\cM}(A) \simeq \HHom_{\HYP^n_R(A)}(C\ten_{R}A,B\ten_RA).
\] 

Of course, if $\fX, \fZ$ are affine schemes, then we can take $n=0$. If they are semi-separated schemes, we can take $n=1$.

For $f:X\ten A \to Z$ in $\cM(A)$,  $\cHom(\fX,\fZ)$  has cohomology groups
\begin{eqnarray*}
\DD^i_{f}(\cHom(\fX,\fZ)/R, M)&=& \EExt^{i+1}_{X\ten A}(\bL^{X\ten A/(X\by Z)\ten A}, \O_{X}\ten_{\H_0R}M)\\
&=&\EExt^{i+1}_{X\ten A}(f^*\bL^{Z/\H_0R}[1], \O_{X}\ten_{\H_0R}M)\\
&=& \EExt^i_Z(\bL^{Z/\H_0R}, \oR f_*\O_{X}\ten_{\H_0R}M).
\end{eqnarray*}

Theorem  \ref{representdaffine} (or a direct application of Theorem \ref{lurierep3} to $\HHom(\fX\by_{\Spec R} \Spec A, \fZ)$)  then shows that   $\cHom(\fX,\fZ)$ is a derived $n$-geometric stack provided  $X$ is strongly quasi-compact (i.e. quasi-compact, quasi-separated \ldots), that $Z$  is almost of finite presentation, and that the groups  $\EExt^i_X(f^*\bL^{Z/\H_0R},  \O_{X}\ten_{\H_0R}-)$ have suitable finiteness properties (which are satisfied if $X$ is proper).
\end{example}

When $\fX$ is a smooth projective scheme and $G=\GL_r$, the following appears as \cite{hag2} Corollary 2.2.6.15:

\begin{example}[Moduli of torsors]\label{Gtorsors}
A special case of Example \ref{Mor} is to take a smooth group scheme $G$ over $R$, and to set $\fZ$ to be the classifying stack $BG$. Then 
$\cHom(\fX, BG)$ is the derived stack of $G$-torsors on $\fX$ --- when $G= \GL_r$, this gives derived moduli of vector bundles, and when $G=\SL_r$ it gives derived moduli of vector bundles with trivial determinant.  

To calculate the cohomology groups $\DD^i_f(\tilde{M}/R,M)$ for $f: \fX \by \Spec A \to BG$,  write $e: \Spec R \to G$ for the unit. Then $e^*\Omega(G/R)$ has a $G$-action, so can be regarded as a quasi-coherent sheaf $\sC$ on $BG$. We then have $\bL^{BG/R}\simeq \sC[1]$, so 
\[
\DD^i_f(\cHom(\fX, BG) /R,M)= \Ext^{i+1}_{X\ten A}(f^*\sC, \O_X\ten_RM).
\]
If $f$ corresponds to a $G$-torsor $P$, then note that $f^*\sC= P\by_G(e^*\Omega(G/R))$.

If $G= \GL(V)$, then $\sC\cong \End(V)$ with the adjoint action, and for $f:\fX \by \Spec A \to B\GL_n$ corresponding to a vector bundle $\sE$, we have $f^*\sC \cong \sEnd(\sE)$, so
\begin{eqnarray*}
\DD^i_f(\cHom(\fX, BG)/R,M)&= &\Ext^{i+1}_{X\ten A}(\sEnd(\sE) , \O_X\ten_RM)\\
&= &\Ext^{i+1}_{X\ten A}(\sE , \sE\ten_RM).
\end{eqnarray*}
\end{example}

\begin{example}[Moduli of polarised schemes]\label{modpolsch}
\label{polar}
For an example combining aspects of all the problems so far, consider moduli of polarised schemes.
In Theorem \ref{representdaffine}, take $\fY= B\bG_m$, and let $\cM(A)$ be some groupoid of flat schemes $X$ over $A$, equipped with polarisations $f: X \to B\bG_m$, with morphisms required to respect the polarisations. Then we have $\tilde{\cM}$ as defined at the the beginning of the section. Again, $\bar{W}\tilde{\cM}$ will be a geometric derived $1$-stack whenever $\cM$ is an algebraic $1$-stack and the cohomology groups $\DD^i(\bar{W}\tilde{\cM}/R, -)$ satisfy the requisite finiteness conditions.

Explicitly,
\[
\DD^i_X(\bar{W}\tilde{\cM}/R, M)= \EExt_X^{i+1}(\bL^{X/B\bG_m},\O_X\ten_{\H_0R}M).
\]
The calculations of the previous example show that  $f^*\bL^{B\bG_m}= \O_X[1]$, so   the cotangent complex $\bL^{X/B\bG_m} $ is  an extension of $\O_X$ by $ \bL^{X/\H_0R} $. This gives a long exact sequence
\[
\ldots \H^{i+1}(X,\O_X\ten_{\H_0R}M) \to\DD^i_X(\bar{W}\tilde{\cM}/R, M) \to \EExt_X^{i+1}(\bL^{X/\H_0R},\O_X\ten_{\H_0R}M)\to \ldots, 
\]
so properness of $X$ gives finiteness.
\end{example}

\begin{remark}
When considering moduli of polarised varieties, it is common only to consider the polarisation up to equivalence, so $(X, \cL^m) \sim (X, \cL^n)$. The corresponding derived moduli stack can be constructed by taking the filtered colimit
\[
\LLim_{n \in \N^{\by}}\bar{W}\tilde{\cM},
\]
for $\tilde{\cM}$ as in Example \ref{polar}, with the transition maps $\tilde{\cM} \to \tilde{\cM}$ (corresponding to $m \to mn$ in $\N^{\by}$) coming from the multiplication map $[n]: \bG_m \to \bG_m$. 

Alternatively, we could simply replace $\bG_m$ with the  (non-quasi-separated) algebraic group space
$\LLim_{n \in \N^{\by}} \bG_m= \bG_m/\mu_{\infty}$, where $\mu_{\infty}= \bigcup_{n \in \N^{\by}}\mu_n$.
\end{remark}

\begin{example}[Moduli of group schemes]\label{modgpsch}
Given a smooth group scheme $G$ over $A$, the nerve $BG$ is a $1$-truncated geometric Artin stack over $A$. However, beware that $\HHom_A(BG, BG')$ is the nerve of the groupoid $[\Hom_{\gp\Sch}(G,G')/G'(A)]$, where $\Hom_{\gp}$ denotes group scheme homomorphisms and $G'$ acts by conjugation. If $\cM(A)$ is the nerve of some open $2$-subcategory of the $2$-category of group schemes over $A$, then $B$ gives a map $\cM \to F\Stack(A)$, so we may apply Theorem \ref{representdaffine} to give a derived moduli functor $\tilde{\cM}$ whenever $\cM$ satisfies the conditions of that theorem.

Since we usually wish to kill the $2$-morphisms (given by inner automorphisms), an alternative approach is to apply Theorem \ref{augrepresentdaffine}, taking $\fZ=\fY= \Spec R$. Making use of the canonical basepoint $i: \Spec A \to BG$, we then have
\[
\HHom_A(BG, BG')\by_{i^*,\HHom_A(\Spec A, BG')} \{i'\} \simeq \Hom_{\gp\Sch}(G,G'),
\]
so the nerve functor allows us to regard smooth group schemes over $A$ as an open subcategory of $d\Stack^{\Spec R}_{\Spec R}(A)$.  

As in Example \ref{Gtorsors},  $\bL^{BG/A}\simeq \sC[1]$, where $\sC= e^*\Omega(G/A)$ with its canonical $G$-action. Thus 
\[
\DD^i_G(\bar{W}\tilde{\cM}/R, M) \cong \Ext^{i+2}_{BG}(\sC, \O_{BG}\ten_A^{\oL}M)
\]
for $i \ge 0$ (or, in the unpointed case, for all $i$). For $i<-1$, the groups  $\DD^i_G(\bar{W}\tilde{\cM}/R, M)$ are $0$ in the pointed case, with $i=-1$ giving an exact sequence
\begin{eqnarray*}
0 \to &\Hom_{BG}(\sC, \O_{BG}\ten_A^{\oL}M) \to \Hom_A(e^*\sC,M)\to \\
&\DD^{-1}_G(\bar{W}\tilde{\cM}/R, M)\to   \Ext^{1}_{BG}(\sC, \O_{BG}\ten_A^{\oL}M)\to 0.
\end{eqnarray*}
\end{example}

\begin{remark}
 Given a flat finitely presented group scheme $G$ over $A$, \cite{Artin} Theorem 6.1 implies that $BG$ is  a $1$-truncated geometric Artin stack over $A$. The analysis of the example above then carries over, subject to replacing $\Omega_G$ with the complex $\chi_G$ defined in   \cite{illusiethesis} \S 2.5.1. The same substitution will also work in the following example, with \cite{toenflat} Theorem 0.1 ensuring that  $K(G,n)$ be  an $n$-truncated geometric Artin stack over $A$.   
\end{remark}

\begin{example}[Moduli of commutative group schemes]\label{commgp}
For smooth commutative group schemes $G$ over $A$, we can no longer study moduli via $BG$, since this has non-commutative deformations. However, for all $n \ge 2$, the construction $K(G,n)$ gives an $n$-truncated geometric Artin stack over $A$, whose deformations all correspond to commutative group schemes. Rather than fixing one value of $n$, the most natural approach is to use them all simultaneously.

The loop space functor $(X,x) \mapsto X^{S^1}\by_X\{x\}$ gives a endofunctor $\Omega$ of $d\Stack_{\Spec R}^{\Spec R}$. This maps $K(G,n)$ to $K(G, n-1)$ (up to weak equivalence). We define $d\Sp\Stack_R(A)$ to be the homotopy inverse limit
of the diagram
\[
\ldots \xra{\Omega} d\Stack_{\Spec R}^{\Spec R} \xra{\Omega} d\Stack_{\Spec R}^{\Spec R}(A);
\] 
this is the $\infty$-category of spectral derived stacks. Defining $F\Sp\Stack_R(A)$ from $F\Stack_R^R(A)$ analogously, we have a functor from the category $S\Gp\Sch(A) $ smooth commutative group schemes $G$ over $A$ to  $F\Sp\Stack_R(A)$, given by $G \mapsto \{K(G,n)\}_n$. Explicitly, $K(G,0)=G$, and we form $K(G,n)$ from the group $n-1$-stack $K(G, n-1)$ as the nerve $K(G,n):= \bar{W}K(G,n-1)$. Then  $S\Gp\Sch$ is an open subcategory of $F\Sp\Stack_R(A)$, each of the projections $K(-,n):S\Gp\Sch\to d\Stack_{\Spec R}^{\Spec R}(A) $ also being open for $n \ge 2$.

We can now adapt Theorem \ref{augrepresentdaffine}, replacing $d\Stack_{\fY}^{\fZ}$ with   $d\Sp\Stack$, noting that the latter is homotopy-preserving and homotopy-homogeneous.
 In order to determine which open subcategories $\cM$  of $S\Gp\Sch $ yield derived geometric stacks
\[
\tilde{\cM}(A):= \cM(\H_0A)\by^h_{F\Sp\Stack(\H_0A)}d\Sp\Stack(A),
\]
we need to calculate the cohomology groups $\DD^i_{G}(\bar{W}\tilde{\cM}/R, M)$. 
The theorem implies that these groups are given by cohomology of
\begin{eqnarray*}
&& \Lim_n\oR\HOM_{K(G,n)}(\bL^{K(G,n)/A}, \O_{K(G,n)}\ten_A^{\oL}M \to \oR i_* M)[1]\\
&=& \Lim_n\oR\HOM_{K(G,n)}(\O_{K(G,n)}\ten_Ae^*\Omega(G/A)[n+1], \O_{K(G,n)}\ten_A^{\oL}M \to \oR i_* M)\\
&=& \oR\HOM_A(e^*\Omega(G/A), \Lim_n\oR \tilde{\Gamma}(K(G,n), \O_{K(G,n)})\ten_A^{\oL}M)[n+1],
\end{eqnarray*}
where $\oR \tilde{\Gamma}= \oR\Gamma/\H^0$. We can rewrite this as
\[
\DD^i_{G}(\bar{W}\tilde{\cM}/R, M) = \EExt^{i+1}_{\Z}( \LLim_n N \Z K(G,n)[n], \Hom_A(e^*\Omega(G/A),M)),
\] 
where $N$ denotes Dold--Kan normalisation, and $\EExt$ is taken in the category of abelian fppf sheaves over $\Spec A$. This ties in with the approach of \cite{Ill2} \S VII.4 (as sketched in \cite{illusiethesis} Theorem 5.3).

Now, \cite{Cartan2} Theorem 9.7 (adapted to abelian sheaves as in  \cite{breeneilmac} Theorem 3) shows that
\[
\LLim_n\sH_{n}( N \Z K(G,n))=G,
\]
and that for fixed $q$, the groups $\sH_{n+q}( N \Z K(G,n))$ are constant for $n>q$, and can be expressed in terms of $G/pG$ and ${}_pG$ for primes $p$. [Explicit descriptions up to q=5 are given at the end of \cite{breeneilmac} \S 1.]
In particular, $\sH_{n+1}(\LLim_n N \Z K(G,n))=0$ for $n>1$, so
\[
\DD^i_{G}(\bar{W}\tilde{\cM}/R, M) = \Ext^{i+1}_{\Z}(G, \Hom_A(e^*\Omega(G/A),M))
\]
for all $i \le 0 $. 
\end{example}

\begin{remark}
This approach means that we can think of (derived) commutative group schemes as functors from simplicial (or dg) algebras to spectra. By contrast, derived stacks are functors from simplicial  algebras to simplicial sets. This might help to explain why Shimura varieties are one of the rare examples of moduli spaces which extend naturally to brave new stacks (i.e. functors from ring spectra to simplicial sets), as observed in \cite{lurieellipticsurvey}.
\end{remark} 

\begin{example}[Moduli of abelian schemes]\label{abvar}
A special case of the previous example is when $\cM$ is an open subcategory of the category of abelian schemes. The description of the cohomology groups in this case simplifies greatly. 
The key observation is that 
\[
\oR \Gamma (G, \O_G) \simeq \Symm_A(\H^1(G,\O_G)[-1]),
\] 
with $\H^1(G,\O_G)$ projective over $A$. 

If we have  $\H^1(G,\O_G)\cong A^g$, then the equivalence becomes    
\[
\oR \Gamma (G, \O_G) \simeq \CC^{\bt}(K(\Z^g,1),\Z)\ten_{\Z}^{\oL}A.
\]
Thus 
\begin{eqnarray*}
\Lim_n\oR \Gamma(K(G,n), \O_{K(G,n)})[n] &\simeq& \Lim_n \CC^{\bt}(K(\Z^g,n+1),\Z)[n]\ten_{\Z}^{\oL}A\\
&=&\CC^{\bt-1}_{\st}(H\Z^g,\Z)\ten_{\Z}^{\oL}A\\
&=& \CC^{\bt-1}_{\st}(H\Z,\Z)\ten_{\Z}^{\oL}A^g,
\end{eqnarray*}
the complex of stable cochains of the integral Eilenberg--Mac Lane spectrum $H\Z$.

In general, $\Lim_n\oR \Gamma(K(G,n), \O_{K(G,n)})[n]$ can be rewritten $A$-linearly in terms of the   cocommutative $A$-coalgebra $\oR \Gamma(G,\O_G) $ (a model for which is the normalised \v Cech complex with Eilenberg-Zilber comultiplication). This gives us  a morphism
\[
\CC^{\bt}_{\st}(H\Z,\Z)\ten_{\Z}^{\oL}\Lim_n\oR \Gamma(K(G,n), \O_{K(G,n)})[n] \to \Lim_n\oR \Gamma(K(G,n), \O_{K(G,n)})[n],
\]
inducing a map
\[
\CC^{\bt-1}_{\st}(H\Z,\Z)\ten_{\Z}^{\oL}\H^1(G,\O_G) \to \Lim_n\oR \Gamma(K(G,n), \O_{K(G,n)})[n].
\]
We know that this is an isomorphism locally on $A$ (as $\H^1(G,\O_G)$ is locally free), so it is an isomorphism in general.

 Thus 
\[
\DD^i_{G}(\bar{W}\tilde{\cM}/R, M)\cong    \Hom_A(e^*\Omega(G/A),\H^1(G,\O_G))\ten_A \H^i_{\st}(H\Z,M).
\]

These stable cohomology groups can be calculated by applying the universal coefficient theorem to   \cite{Cartan2} Theorem 7; alternatively we could apply the  universal coefficient theorem for $\Tor$ to the  integral Steenrod algebra $\H^i_{\st}(H\Z,\Z)$. The important features are that $\H^{<0}_{\st}(H\Z,M)=0$, while $\H^0_{\st}(H\Z,M)=M$ and $\H^{>0}_{\st}(H\Z,M)$ is an invariant of the torsion and cotorsion in $M$. In particular, if $M$ is a $\Q$-module, then $\H^{>0}_{\st}(H\Z,M)=0$. We always have $\H^{1}_{\st}(H\Z,M)=0$, which is why moduli of abelian varieties are unobstructed.
\end{example}

Beware that the example above only gives a derived Artin stack when a suitable moduli stack of abelian varieties exists. 
\begin{remark}\label{ppav}
 Most approaches to moduli of abelian varieties involve polarisations, and we now sketch a possible approach. A (principal) polarisation is a symmetric bilinear map $G \by G \to \bG_m$, such that the restriction to the diagonal $G \to G \by G$ gives an ample line bundle on $G$, and the the induced map $G \to G^{\vee}$ is an isogeny (resp. an isomorphism).

For an associated derived moduli problem, we could take spectral derived geometric stacks $\fX$ over $A$ equipped with  morphisms
\[
\alpha: [\fX \wedge \fX/S_2] \to \Omega^{-1}H(\bG_m\ten A), 
\]
 where $\Omega^{-1}H\pi= \{K(\pi, n+1)\}_n$, $\wedge$ is the smash product on spectra, and $S_2$ acts by transposition.
 
Since isogenies and isomorphisms are open in the space of group scheme homomorphisms, we could then take $\cM(A)$ to be the category of polarised (or principally polarised) abelian schemes over $A$. For $G \in \cM(A)$,  the associated spectral derived geometric stack $\fX$ would then be  $HG$, with $\alpha$ given by the composition
\[
 [\fX \wedge \fX/S_2] \to [H(G\ten_{\Z}G)/S_2] \to \Omega^{-1}H(\bG_m\ten A),
\]
the right-hand map coming from the polarisation.

In order to verify that this defines a geometric derived stack, the difficult part would be to calculate the cohomology groups in order to check that they satisfy the requisite finiteness conditions. 
\end{remark}

\section{Derived moduli of quasi-coherent sheaves}\label{modshfsn} 

\subsection{Derived quasi-coherent sheaves}

Given a simplicial ring $A_{\bt}$, the simplicial normalisation functor $N^s$ of Definition \ref{normdef} induces an equivalence  between the categories $s\Mod(A)$ of simplicial $A_{\bt}$-modules, and $dg_+\Mod(N^sA)$ of $N^sA$-modules in non-negatively graded chain complexes, where the Eilenberg--Zilber shuffle product makes $N^sA$ into a  graded-commutative algebra. As observed in \cite{hag2} \S 2.2.1, this extends to give a weak equivalence between the $\infty$-category of stable $A$-modules, and the $\infty$-category $dg\Mod(N^sA)$ of $N^sA$-modules in $\Z$-graded chain complexes, and hence an equivalence between the corresponding homotopy categories. More precisely, the $\infty$-structure on $dg\Mod(N^sA)$ comes from a model structure:

\begin{definition}
  In  $dg\Mod(N^sA)$, fibrations are surjections, weak equivalences are quasi-isomorphisms, and cofibrations satisfy LLP with respect to trivial fibrations.   Explicitly, $f: L \to M$  is a a cofibration if it is a retract of a transfinite composition of monomorphisms $f':L' \to M'$ with $\coker f' \cong (N^sA)[n]$ for some $n$.
\end{definition}
In particular, if $A$ is discrete and $f$ is a monomorphism with $\coker f$ bounded below and levelwise projective, then $f$ is a cofibration.

\begin{definition}
 Given $L, M \in dg\Mod(N^sA)$, we can then define $\HOM_A(L,M)$ to be the chain complex given in level $n$ by $A$-linear maps $L \to M[n]$ of graded objects. The differential on $ \HOM_A(L,M)$ is given by $d(f)= d\circ f \pm f\circ d$. The simplicial structure on $dg\Mod(N^sA) $ is then given by setting
\[
\HHom_A(L,M):= (N^s)^{-1}(\tau_{\ge 0}\HOM_A(L,M)),
\]
where $\tau_{\ge 0}$ is good truncation.
\end{definition}
The associated $\infty$-category is given by the full simplicial subcategory of $ dg\Mod(N^sA)$ on cofibrant objects. Equivalently, we have derived function complexes on $ dg\Mod(N^sA)$ given by  $\oR\HHom(L,M):= \HHom(\tilde{L},M)$, for $\tilde{L} \to L$ a cofibrant replacement.

\begin{definition}
An object $X_{\bt}$ in $sd\Aff_R$ is said to be a  homotopy derived Artin $n$-hypergroupoid if it admits a map $X_{\bt} \to \tilde{X}_{\bt} $ to a   derived Artin $n$-hypergroupoid with the property that the level maps $X_n \to   \tilde{X}_n$ are all weak equivalences in $d\Aff_R$.
\end{definition}

\begin{definition}
If  $X$ is a  homotopy derived Artin $l$-hypergroupoid, define  the simplicial cosimplicial algebra  $O(X)^{\bt}_{\bt}$ by $O(X)^n_m= \Gamma(X_n^m, \O_{X_n^m})$. Define the category $dg\Mod(X)$ to have objects $M$ 
 consisting of complexes $M^n_{\bt}\in  dg\Mod(N^sO(X)^n_{\bt})$ for all $n$, together with morphisms
\begin{eqnarray*}
\pd^i:M^n_{\bt}\ten_{N^sO(X)^n_{\bt}, \pd_i^*}N^sO(X)^{n+1}_{\bt} &\to& M^{n+1}_{\bt}\\
\sigma^i: M^n_{\bt}\ten_{N^sO(X)^n_{\bt}, \sigma_i^*}N^sO(X)^{n-1}_{\bt} &\to& M^{n-1}_{\bt},
\end{eqnarray*}
satisfying the usual cosimplicial identities.
\end{definition}

\begin{definition}
A morphism $f: M \to N$ in $dg\Mod(X)$ is said to be a weak equivalence if each map $f^n: M^n \to N^n$ is a weak equivalence in  $dg\Mod(N^sO(X)^n_{\bt})$. Say that $f$ is a fibration if each map  $f^n: M^n \to N^n$ is surjective.

Finally, $f$ is  said to be a cofibration if it satisfies the LLP with respect to trivial fibrations. For an explicit characterisation, first form  the $n$th latching object $L^nM \in dg\Mod(N^sO(X)^n_{\bt})$ as the cokernel of 
\[
(\alpha -\beta): \bigoplus_{i=0}^{n} \bigoplus_{j=0}^{i-1} \pd_i^*\pd_j^*M^{n-2} \to \bigoplus_{i=0}^{n-1} \pd_i^*M^{n-1},
\]
where for $x \in \pd_i^*\pd_j^*M^{n-2}= \pd_j^*\pd_{i-1}^*M^{n-2}$, we define $\alpha(x)= \pd_i^*(\pd^j)x \in \pd_i^*M^{n-1}$, and $\beta(x)= \pd_j^*(\pd^{i-1})x \in \pd_j^*M^{n-1}$. Then $f$ is a cofibration if and only if the maps
\[
L^nN\oplus_{L^mM}M^n \to N^n
\]
are cofibrations in $ dg\Mod(N^sO(X)^n_{\bt})$ for all $n$.
\end{definition}

We can then form the simplicial structure by analogy with $dg\Mod(N^sA)$: define the chain complex
\[
\HOM_X(M,N)_n:= \Hom_{O(X)}(M, N[n]),
\]
where $\Hom$ respects cosimplicial, but not chain, structures. Then $\HHom(M,N):= (N^s)^{-1}(\tau_{\ge 0} \HOM_X(M,N))$, with $\oR \HHom$ given by taking a cofibrant replacement for $M$. 

\begin{definition}
Define $dg\Mod(X)_{\cart} \subset dg\Mod(X)$  to be the full subcategory consisting of those  $M$ for which the morphisms $\pd^i$ are all weak equivalences. Let  $\Ho(dg\Mod(X))_{\cart})$ be the category obtained by localising at weak equivalences. Note that $dg\Mod(X)_{\cart}$ also inherits a simplicial structure from $dg\Mod(X) $.
Objects of $dg\Mod(X)_{\cart}$ are called  derived quasi-coherent sheaves  on $X$. 

Also define $dg\Mod(X)_{\cart,>-\infty}\subset dg\Mod(X)_{\cart}$ to be the full subcategory consisting of objects $M$ for which the chain complexes $M^n$ are all bounded below (i.e. bounded above if regarded as cochain complexes). 
\end{definition}

\begin{proposition}\label{dsheafequiv}
Fix a  homotopy derived Artin $n$-hypergroupoid $X$.
The  simplicial category of derived quasi-coherent sheaves on $X$ is equivalent to the $\infty$-category $dg\Mod_{\cart}(X^{\sharp})$  of  homotopy-Cartesian modules on $X^{\sharp}$ (as in \cite{hag2} Definition 1.2.12.1) or equivalently, of  quasi-coherent complexes on $X^{\sharp}$ (in the sense of \cite{lurie} \S 5.2).
\end{proposition}
\begin{proof}
This is \cite{stacks2} Corollary \ref{stacks-qcohequiv}.
\end{proof}
In particular, this means that $dg\Mod(X)_{\cart} $ is essentially independent of the atlas $X$ chosen for $X^{\sharp}$.

\begin{definition}\label{npdef}
Define an  Artin  (resp. Deligne--Mumford) $n$-hypergroupoid  over a ring $R$ to be an object  $X_{\bt}\in s\Aff_R$, such that the    partial matching maps
$$
X_m \to M_{\L^m_k} X 
$$
are   smooth (resp. \'etale)  covers  for all $k,m$, and isomorphisms  for all $m>n$ and all $k$.
\end{definition}
In particular, note that every $n$-hypergroupoid is a  homotopy derived  $n$-hypergroupoid. Specialising to underived $n$-stacks, Proposition \ref{bigthmd} shows that every strongly quasi-compact $n$-geometric Artin  (resp. Deligne--Mumford) stack is of the form $X^{\sharp}$, for some  Artin  (resp. Deligne--Mumford) $n$-hypergroupoid $X$.

For our purposes, the main applications of derived quasi-coherent sheaves stem from the following:
\begin{proposition}\label{cringmod}
\begin{enumerate}
\item Take an  Artin $n$-hypergroupoid $X$. Then the category of quasi-coherent sheaves on the geometric stack $X^{\sharp}$ is equivalent to the full subcategory of $ \Ho(dg\Mod(X)_{\cart} )$  consisting of complexes concentrated in chain degree $0$.

\item If $Y$ is a semi-separated quasi-compact scheme, take an open affine cover $\{U_i\}_{i=1}^n$. Let $U= \coprod_i U_i$,
and set $X$ to be the  Deligne-Mumford $1$-hypergroupoid given by
\[
X_n = \overbrace{U\by_Y U \by_Y \ldots \by_YU}^{n+1}.
\]
Then  $ \Ho(dg\Mod(X)_{\cart} )$ is equivalent to the derived  category of complexes of quasi-coherent sheaves on $Y$.
\end{enumerate}
\end{proposition}
\begin{proof}
The first statement is \cite{stacks2} Corollary \ref{stacks-qcohequiv}. Once we note that $X^{\sharp}=Y$, the second statement follows from the observation in the introduction of \cite{huettemann} that the  proof of  \cite{huettemann} Theorem 4.5.1 works for any   quasi-compact
scheme equipped with a finite semi-separating affine covering.
\end{proof}

\subsection{Moduli}

\begin{definition}
Take  $A \in dg_+\Alg_R$ (resp. $A \in s\Alg_R$), and a  homotopy derived Artin $n$-hypergroupoid $X$ over $R$, chosen so that each $X_n$ is levelwise flat over $R$ (which holds for instance if $X$ is a derived Artin $n$-hypergroupoid). Then observe that $X\ten_RA$ is a  homotopy derived Artin $n$-hypergroupoid over $A$, and
define
\[
d\CART_X(A), \quad \text{ resp. }\quad  d\CART_{X,>-\infty}(A),
\]
to be the full simplicial subcategory of the category  $dg\Mod(X\ten_RA)_{\cart}$,  resp. $dg\Mod(X\ten_RA)_{\cart,>-\infty}$  on cofibrant  objects. 
\end{definition}

\begin{proposition}\label{cmodhgs}
The functor $d\CART_{X,>-\infty}:  d\cN^{\flat} \to s\Cat$ defined above is 2-homogeneous and formally 2-quasi-smooth.
\end{proposition}
\begin{proof}
The proof of Proposition \ref{calghgs} carries over to this context, with cofibrancy once again ensuring the existence of suitable lifts; the bounded below hypothesis ensures that a flat deformation of a cofibrant object is again cofibrant, since it is quasi-free and we can order generators by degree. 
\end{proof}

\begin{theorem}\label{representdmod}
Take a strongly quasi-compact $m$-geometric  derived Artin stack $\fY$   over $R$, and choose a  homotopy derived Artin $n$-hypergroupoid $\check{\fY}$ with $\check{\fY}^{\sharp}=\fY$ and each $\check{\fY}_n$ flat over $R$.

Assume that we have  an  $n$-truncated presheaf $\cM: \Alg_{\H_0R} \to s\Cat$, open in the functor $ A \mapsto \cW(dg\Mod_{\cart,>-\infty}(\fY\ten_R^{\oL}A))$. 

Also assume that this satisfies the following conditions:
\begin{enumerate}
\item[(0)] if $\{f_i \co \Spec B_i \to \Spec A\}_{i \in I}$ is an \'etale cover in affine $\H_0R$-schemes, then an object   $\sE \in \pi_0\cW(dg\Mod_{\cart}(\fY\ten_R^{\oL}A))$ lies in the essential image of $\pi_0\cM(A)$ whenever $f_i^*\sE$ lies in the essential image of $\pi_0\cM(B_i)$ for all $i$.

\item\label{dmafp2} for all finitely generated $A \in \Alg_{\H_0R}$  and all $\sE \in \cM(A)$, the functors 
\[
 \EExt^i_{\fY\ten_R^{\oL}A}(\sE,\sE\ten_A^{\oL}- ) :\Mod_A \to \Ab
\]
 preserve filtered colimits for all $i\ne 1$.

\item for all finitely generated integral domains $A \in \Alg_{\H_0R}$  and all $\sE \in \cM(A)$, the groups $\EExt^i_{\fY\ten_R^{\oL}A}(\sE, \sE )$ are all  finitely generated $A$-modules.

\item
The functor $c(\pi_0\cM):  \Alg_{\pi_0R} \to \Set$  of components of the groupoid $\pi_0\cM$ preserves filtered colimits.

\item\label{Ahatcdn} for all complete discrete local Noetherian  $\H_0R$-algebras $A$, with maximal ideal $\m$, the map
\[
c(\pi_0\cM(A)) \to \Lim_r  c(\pi_0\cM(A/\m^r))
\] 
is an isomorphism, as are the maps 
\[
 \EExt^i_{\fY\ten_R^{\oL}A}(\sE, \sE ) \to \Lim_r \EExt^i_{\fY\ten_R^{\oL}A}(\sE, \sE/\m^r )
\]
 for  all $\sE \in \cM(A)$ and all $i\le 0$.

\end{enumerate}

Let $\tilde{\cM}:d\cN^{\flat}_R \to s\Cat$  be the full simplicial subcategory of $\cW(d\CART_{\check{\fY} }(A))$ consisting of objects $\sF$ for which $\sF\ten_A\H_0A$ is weakly equivalent  in $dg\Mod_{\cart}(\fY\ten_R^{\oL}\H_0A)$  to an object of  $\cM(\H_0A)$. 
Then the nerve $\bar{W}\cM$ is an geometric derived $n$-stack.
\end{theorem}
\begin{proof}
First note that Proposition \ref{dsheafequiv} allows us to replace $d\CART_{\check{\fY} }(A)$ with $ dg\Mod_{\cart}(\fY\ten_R^{\oL}A) $, so
\[
\tilde{\cM}(A) \simeq \cM(\H_0A)\by^h_{\cW(dg\Mod_{\cart}(\fY\ten_R^{\oL}\H_0A) )}\cW(dg\Mod_{\cart}(\fY\ten_R^{\oL}A) ),
\] 
because $\cM $ is open in $ \pi^0\cW(dg\Mod_{\cart}(\fY\ten_R^{\oL}-))$. 

Now, we have
\[
\DD^i_{\sE}(\bar{W}\tilde{\cM}/R, M) =\EExt^{i+1}_{\fY\ten_R^{\oL}A}(\sE, \sE\ten_A^{\oL}M),
\]
 for instance from the calculations of \cite{stacks2} \S \ref{stacks-dsheaves},  where $\EExt$ groups are also characterised as
\[
\EExt^i_{\fY\ten_R^{\oL}A}(\sE, \F) \cong  \H_{-i}\HOM_{X\ten_RA}(\sE,\sF),
\]
provided $\sE$ is chosen to be cofibrant.

Note that $\pi_i\HHom_{\cM(B)}(\sE,\sE)= \EExt^{-i}(\sE, \sE\ten_A^{\oL}B)$ for $i \ge 1$, and that $\pi_0\HHom_{\cM(B)}(\sE,\sE)$ consists of invertible elements in $ \EExt^{0}(\sE, \sE\ten_A^{\oL}B)$.
Thus Condition (\ref{dmafp2}) ensures that
for all $A \in \Alg_{\pi_0R}$ and all $\sE \in \cM(A)$, the functors $\pi_i\HHom_{\cM}(\sE,\sE): \Alg_A \to \Set$  preserve filtered colimits for all $i\ge 0$. 

Next, we need to establish that for all complete discrete local Noetherian  $\H_0R$-algebras $A$, with maximal ideal $\m$, the map
$$
\cM(A) \to {\Lim}^h \cM(A/\m^n)
$$
is a weak equivalence.

Observe that 
\[
 \Aut_{\pi_0\cM(A)}(\sE)=\{(g,h) \in \EExt^0_{\fY\ten_R^{\oL}A}(\sE, \sE )\by  \EExt^0_{\fY\ten_R^{\oL}A}(\sE, \sE ) \,:\, gh=hg=\id\}, 
\] 
so Condition \ref{Ahatcdn} implies that $\Aut_{\pi_0\cM(A)}(\sE) \cong \Lim_r \Aut_{\pi_0\cM(A/\m^r)}(\sE/\m^r) $.
Since
\begin{eqnarray*}
\pi_1(\bar{W}\cM(A),\sE))&=&  \Aut_{\pi_0\cM(A)}(\sE)\\
\pi_i(\bar{W}\cM(A),\sE))&=&  \EExt^{1-i}_{\fY\ten_R^{\oL}A}(\sE, \sE ) 
\end{eqnarray*}
for $i \ge 2$, by  Remark \ref{formalexistrk} it suffices to show that
each system 
\[
\{\im( \pi_{i}(\bar{W}\cM(A/\m^s),x=\sE)\to  \pi_{i}(\bar{W}\cM(A/\m^r),\sE))\}_{s \ge r}
\]
satisfies the Mittag-Leffler condition for $i \ge 1$  and $\sE \in \cM(A)$, and that  Condition (\ref{Ahatcdn}) holds.

Since $\EExt^{1-i}_{\fY\ten_R^{\oL}A}(\sE, \sE/\m^r)$ is a finitely generated $A/\m^r$-module, it is Artinian, so the Mittag-Leffler condition is automatically satisfied for $i \ge 2$. For $i=1$, we can study the obstruction maps of Proposition \ref{obs} to deduce that
\[
\Aut_{\pi_0\cM(A/\m^r)}(\sE)= \EExt^0_{\fY\ten_R^{\oL}A}(\sE, \sE/\m^r )\by_{\EExt^0_{\fY\ten_R^{\oL}A}(\sE, \sE/\m )}\Aut_{\pi_0\cM(A/\m)}(\sE/\m).
\]
Since $\EExt^0_{\fY\ten_R^{\oL}A}(\sE, \sE/\m^r )$ is Artinian, we deduce that the Mittag-Leffler condition holds for $i=1$, so we have satisfied the conditions of Remark \ref{formalexistrk}.

The proof of Corollary \ref{representaffine} now  carries over to this context, using Proposition \ref{cmodhgs} in place of Proposition \ref{alghgs}. For the hypersheaf condition, observe that  the argument of Corollary \ref{representaffine} shows that $\bar{W}\cW$ applied to the full simplicial subcategory of $dg\Mod(X\ten_R-)$ on cofibrant objects gives an \'etale hypersheaf, so $\bar{W}\cW d\CART_X$ and hence $\bar{W}\cW\cM$ are hypersheaves by Proposition \ref{descentlemma}.
\end{proof}

\begin{remarks}
 If $R=\H_0R$, we can let  $\fY=Y$ be a  semi-separated quasi-compact scheme, then  form the \v Cech nerve $\check{Y}$ as in part 2 of Proposition \ref{cringmod}, giving an explicit description of $\tilde{\cM}$. 

In order to satisfy the finiteness hypotheses for $\EExt^*$, a likely candidate for a derived geometric $n$-stack   would be to let $\cM(A)$ consist of perfect complexes $\sE_{\bt}$ of $\O_Y\ten_R^{\oL}A$-modules of  $\Tor$-dimension $(n-1)$ over $A$, with $\sH_i(\sE)$ quasi-coherent of  proper support on $Y$.  
\end{remarks}

\begin{definition}
Given  an  $m$-geometric  (underived) Artin stack $Y$ over $\H_0R$, and $A \in \Alg_{\H_0R}$, let $F\Mod_{Y}(A)$ be the category of quasi-coherent sheaves (not complexes) on $Y\ten_{\H_0R}A$, flat over $A$. 
\end{definition}

\begin{example}\label{modqcoheg}
For $\fY$ as in Theorem \ref{representdmod}, set $Y:= \pi^0\fY$ (an   $m$-geometric  (underived) Artin stack over $\H_0R$).
Of course, if $R=\H_0R$, then we will have $Y=\fY$.
 Choose a  stack $\cM$ equipped with functorial fully faithful embeddings   of $\cM$ into the core of  $F\Mod_{Y}$, with the property that $\cM$ is closed under infinitesimal deformations. This ensures that $\bar{W}\cM \to   \bar{W}F\Mod_{Y}$ is open.

Note that flatness ensures that for any morphism $A \to B$ in $\Alg_{\H_0R}$ and $\sF \in F\Mod_{Y}(A)$, the complex $\sF\ten_A^{\oL}B$ is weakly equivalent to $\sF\ten_AB$, so $F\Mod_Y(A)$ embeds into $\Mod_{\cart}(\fY\ten_R^{\oL}A) )$. In fact, we can characterise $F\Mod_Y(A)$ as being weakly equivalent to the subcategory of $\Mod_{\cart}(\fY\ten_R^{\oL}A)$ consisting of complexes $\sF$ for which the homology sheaves $\sH_i(\sF\ten_A^{\oL}B)=0$ for all $i\ne 0$ and all $B \in \Alg_A$. Thus $F\Mod \to \pi^0\Mod_{\cart}(\fY\ten_R^{\oL}-)$ is an open embedding.

Provided that the objects of $\cM$ satisfy the conditions of Theorem \ref{representdmod}, this gives us a geometric derived $1$-stack of quasi-coherent sheaves of $\fY$.  An example for which the conditions are satisfied arises when $Y$ is a projective scheme over $R$ and $\cM(A)$ consists of coherent sheaves on $Y\ten_R A$, flat over $A$.
\end{example}

\subsubsection{Comparisons}

If we wish to study derived moduli of rank $r$ vector bundles on a  homotopy-flat derived geometric $\infty$-stack $\fX$ over $R$, we now have two possible approaches, which we now show are consistent. Let $X=\fX\ten_R\H_0R$, and define $\cM: \Alg_{\H_0R}\to \gpd$ by letting $\cM(A)$ be the groupoid of rank $r$ vector bundles on $X\ten_{\H_0R}A$. 

Example \ref{Gtorsors} gives one extension $\cHom(\fX, B\GL_r): d\cN^{\flat}_R \to \bS$ of $B\cM$. Theorem \ref{representdmod} gives another candidate  $\bar{W}\tilde{\cM}$ for an extension, where the objects of $\tilde{\cM}:d\cN^{\flat}_R \to s\Cat$   are quasi-coherent complexes $\sE$ on $\fX\ten_R^{\oL}A$ for which $\sE\ten_A^{\oL}\H_0A$ is quasi-isomorphic to a rank $r$ vector bundle on $\fX\ten_R^{\oL}\H_0A$.

The following lemma is essentially the same as \cite{hag2} Proposition 1.3.7.10 and its subsequent remarks.
\begin{lemma}\label{cfbundles}
The derived geometric $1$-stacks $\cHom(\fX, B\GL_r)$ and $\bar{W}\tilde{\cM}$ defined above are weakly equivalent.
\end{lemma}
\begin{proof}
The argument of Example \ref{Gtorsors} gives a morphism from
$\cHom(\fX, B\GL_n)$ to $\bar{W}\cW(\cT)$,  where $\cT(A)$ is the full simplicial subcategory of $d\Stack_{\fX \by B\GL_r}(A)$ on objects $\fZ$ for which the projection $f_1:\fZ \to \fX\ten_R^{\oL}A$ is a weak equivalence. We can then define a morphism $\cW(\cT) \to \tilde{\cM}$ as follows: let  $\sU$ be the universal rank $r$ vector bundle on $B\GL_r$, and send 
\[
( \fZ \xra{(f_1,f_2)} (\fX \by B\GL_r) \ten_R^{\oL}A )
\]
in $\cT(A)$
 to the quasi-coherent complex $\oR f_{1*} f^*\sU$ on $\fX$. 

Thus we have constructed a map $\cHom(\fX, B\GL_n) \to \bar{W}\tilde{\cM}$. It is immediate that this gives an equivalence on $\pi^0$, and straightforward to check that it gives isomorphisms on cohomology groups $\DD^i$, so it is a weak equivalence, by \cite{drep} Proposition \ref{drep-detectweak}. 
\end{proof}

\begin{remark}
We could adapt Theorem \ref{representdmod}  to moduli of polarised varieties, by considering pairs $(T,M)$, for $T \in c\ALG(A)$ and $M \in d\CART_{T}(A)$. The resulting derived stack would parametrise pairs $(\fX,\sE)$, for $\fX$ a derived geometric $m$-stack and $\sE$ a quasi-coherent complex on  $\fX$ (bounded above as a cochain complex). 
The argument of Lemma \ref{cfbundles} then adapts to  show that an open substack of the resulting moduli stack is equivalent to the derived stack from Example \ref{polar}.
\end{remark}

\bibliographystyle{alphanum}
\bibliography{references}
\end{document}